\documentclass[10pt]{article}

\usepackage{amsthm}
\usepackage{amsmath}
\usepackage{amsfonts}
\usepackage{amssymb}
\usepackage{mathrsfs}
\usepackage{empheq}
\usepackage{stmaryrd}
\usepackage{dsfont}
\usepackage{enumerate}
\usepackage{cancel}
\usepackage{comment}
\usepackage{pgfplots}
\usepackage{soul}
\usepackage{multirow}
\usepackage[applemac]{inputenc}
\usepackage{hyperref}
\usepackage{xcolor}

\definecolor{mygreen}{rgb}{0, 0.5, 0}

\newcommand{\1}{\mathds{1}}

\newcommand{\N}{\mathbb{N}}

\newcommand{\R}{\mathbb{R}}
\newcommand{\C}{\mathbb{C}}

\newcommand{\vsS}{\mathbb{S}}

\newcommand{\vH}{{\cal H}}

\newcommand{\Br}{\mathscr{B}}
\newcommand{\Dr}{\mathscr{D}}

\newcommand{\Lr}{\mathscr{L}}

\newcommand{\vphi}{\varphi}
\newcommand{\eps}{\varepsilon}

\newcommand{\dsp}{\displaystyle}
\newcommand{\ovl}{\overline}

\newcommand{\vlim}{\lim\limits}

\newcommand{\vsup}{\sup\limits}

\newcommand{\vint}{\int\limits}

\newcommand{\inj}{\hookrightarrow}
\newcommand{\tends}{\longrightarrow}
\newcommand{\weak}{\rightharpoonup}
\newcommand{\wt}{\widetilde}

\newcommand{\loc}{\mathrm{loc}}

\renewcommand{\b}{\mathrm{b}}
\newcommand{\co}{\mathrm{c}}
\renewcommand{\d}{\mathrm{d}}

\newcommand{\dist}{\mathrm{dist}}

\newcommand{\GN}{\mathrm{GN}}
\newcommand{\vi}{\mathrm{i}}
\newcommand{\w}{{\textsl w}}
\renewcommand{\le}{\leqslant}
\renewcommand{\ge}{\geqslant}
\renewcommand{\Re}{\mathrm{Re}}
\renewcommand{\Im}{\mathrm{Im}}

\newcommand{\p}{\prime}

\newcommand{\eqdef}{\stackrel{\mathrm{def}}{=}}
\DeclareMathOperator{\supp}{supp}

\DeclareMathOperator{\sign}{sign}

\numberwithin{equation}{section}

\newtheorem{thm}{Theorem}[section]
\newtheorem{prop}[thm]{Proposition}
\newtheorem{cor}[thm]{Corollary}
\newtheorem{lem}[thm]{Lemma}

\theoremstyle{definition}
\newtheorem{rmk}[thm]{Remark}
\newtheorem{defi}[thm]{Definition}

\newtheorem{exa}[thm]{Example}
\newtheorem{assum}[thm]{Assumption}

\newenvironment{proof*}{\noindent{\bf Proof.}}{\qed}
\newenvironment{vproof}[1]{\noindent{\bf Proof #1}}{\qed}

\title{\huge Strong stabilization of damped nonlinear Schrödinger equation with saturation on unbounded domains}
\author{\sc Pascal Bégout$^*$ and Jes\'us Ildefonso D{\'{\i}}az$^\dagger$}
\date{}

\begin{document}

\maketitle

\begin{gather*}
\begin{array}{cc}
	^*\text{Toulouse School of Economics}	&	\;^\dagger\text{Instituto de Matem\'atica Interdisciplinar}	\\
		\text{Université Toulouse Capitole}	&	\text{Universidad Complutense de Madrid}			\\
\text{Institut de Mathématiques de Toulouse}	&	\text{Plaza de las Ciencias, 3}						\\
		\text{1, Esplanade de l’Université}	&	\text{28040 Madrid, SPAIN}						\\
	\text{31080 Toulouse Cedex 6, FRANCE}	&
\bigskip \\
\text{
{\footnotesize E-mail\:: \href{mailto:Pascal.Begout@math.cnrs.fr}{\texttt{Pascal.Begout@math.cnrs.fr}}}}
&
\text{
{\footnotesize E-mail\:: \href{mailto:jidiaz@ucm.es}{\texttt{jidiaz@ucm.es}}}
}
\end{array}
\end{gather*}

\begin{abstract}
We consider the damped nonlinear Schr\"{o}dinger equation with saturation: i.e., the complex evolution equation contains in its left hand side, besides the potential term $V(x)u,$ a nonlinear term of the form $\vi\mu u(t,x)/|u(t,x)|$ for a given parameter $\mu >0$ (arising in optical applications on non-Kerr-like fibers). In the right hand side we assume a given forcing term $f(t,x).$ The important new difficulty, in contrast to previous results in the literature, comes from the fact that the spatial domain is assumed to be unbounded. We start by proving the existence and uniqueness of weak and strong solutions according the regularity of the data of the problem. The existence of solutions with a lower regularity is also obtained by working with a sequence of spaces verifying the Radon-Nikod\'{y}m property. Concerning the asymptotic behavior for large times we prove a strong stabilization result. For instance, in the one dimensional case we prove that there is extinction in finite time of the solutions under the mere assumption that the $L^\infty$-norm of the forcing term $f(t,x)$ becomes less than $\mu$ after a finite time. This presents some analogies with the so called feedback \textit{bang-bang controls} $v$ (here $v=-\mathrm{i}\mu u/|u|+f).$
\end{abstract}

{\let\thefootnote\relax}
{\let\thefootnote\relax}

\tableofcontents

\baselineskip .6cm

\section{Introduction}
\label{introduction}

The main goal of this paper is the consideration of the following damped nonlinear Schrödinger equation
\begin{empheq}[left=\empheqlbrace]{align}
	\label{nls}
	\vi\frac{\partial u}{\partial t}+\Delta u+V(x)u+\vi\mu\frac{u}{|u|}=f(t,x),		&	\text{ in } (0,\infty)\times\Omega,				\\
	\label{nlsb}
	u_{|\partial\Omega}=0,										&	\text{ on } (0,\infty)\times\partial\Omega,	\dfrac{}{}	\\
	\label{u0}
	u(0)= u_0,													&	\text{ in } \Omega,
\end{empheq}
where $\vi^2=-1,$ $\mu>0,$ $\Omega\subseteq\R^N,$ $f\in L^1_\loc\big([0,\infty);L^2(\Omega)\big),$ $V\in L^1_\loc(\Omega;\R)$ and $u_0\in L^2(\Omega).$

It is important to point out that although the study of some damped nonlinear Schrödinger equations $(\mu>0)$ was already considered since the seventies of the past century, in most of the cases the domain $\Omega$ was assumed to be an open bounded subset of $\R^N,$ or the nonlinear term was usually assumed to be a Lipschitz continuous function (see Cazenave~\cite{MR2002047} and the references therein, and references, for instance, in \cite{MR4053613}). We will explain later why the case of $\Omega$ unbounded presents important difficulties in its treatment.

The nonlinear term considered in this paper (see \eqref{nls}) is usually called as a ``saturation nonlinearity". This kind of nonlinearity arises quite often in the modeling of some problems in the framework of optical applications in non-Kerr-like fibers (see, e.g., Gatz and Herrmann~\cite{gh}, Lyra and Gouveia-Neto~\cite{lg}, Tatsing, Mohamadou and Tiofack~\cite{TMT}, and the references given in \cite{MR4053613}).

The case of a saturation nonlinearity as the one considered in \eqref{nls} can be understood also in the framework of Control Theory as a special case of feed-back control of ``bang-bang type"
\begin{gather}
\label{y}
y(t,x)=\vi\mu\frac{u(t,x)}{|u(t,x)|}.
\end{gather}
This type of control has been considered in the applications to many dissipative evolution equations (see, e.g. Lasiecka and Seidman~\cite{MR2020641}, Tarbouriech, Garcia, Gomes da Silva and Queinnec~\cite{MR3024786}, and Laabissi and Taboye~\cite{MR4303962}). Nevertheless, the controllability for Schrödinger equations is more delicate (see, e.g., Machtyngier~\cite{MR1255957}) and requires ``ad hoc" arguments.

The main goal of this paper is to offer a mathematical treatment of problem \eqref{nls} for the case that was left as an open problem in some previous results in the literature (see, e.g., Carles and Gallo~\cite{MR2765425}, Carles and Ozawa~\cite{MR3306821}, Hayashi~\cite{MR3802567} and \cite{MR4053613,MR4340780,MR4503241}): the case of $\Omega$ an unbounded domain. More precisely, the case $m=0$ is not treated in Hayashi~\cite{MR3802567} and in \cite{MR4053613,MR4503241}. The case $m=0$ and $\Omega$ unbounded (actually $\Omega=\R^N)$ is partially treated in Carles and Ozawa~\cite{MR3306821}: $m\in[0,1]$ and $N=1.$

Before describing the main difficulties in the study of this case, let us point out that the relevance of the consideration of this formulation comes from the fact that we will prove (in Section~\ref{proofext}) that thanks to this kind of nonlinearity (i.e., for controls $y\in L^\infty\big((0,\infty)\times\Omega;\C\big)$ given by \eqref{y}) we will show the stabilization to zero (as $t\to\infty)$ of solutions of \eqref{nls} even if the source term $f(t,x)$ is merely a bounded function (i.e., without require that $f(t,\:.\:)\to0,$ in some sense, as $t\to\infty).$ Moreover, in the case of $N=1,$ we will prove that the finite time extinction (which under suitable conditions on $f(t,\,.\,))$ may be instantaneous -- see Theorem~\ref{thmextN1}, part~\ref{thmextN14}.

A natural way to start the study of equation \eqref{nls} is to enlarge the framework by considering the more general sublinear Schrödinger equation (now not necessarily saturated) of the form
\begin{gather}
\label{nlsg}
\vi\frac{\partial u}{\partial t}+\Delta u+V(x)u+a|u|^{-(1-m)}u=f(t,x), \text{ in } (0,\infty)\times\Omega,
\end{gather}
where $a\in\C$ and $m\in[0,1].$ Such as we will indicate later, different authors introduced some constraints among the parameters $m$ and $a$ in order to get the existence and uniqueness of different types of solution (this will be properly presented in the Section~\ref{exiuni} below). So, if $m\in[0,1]$ then we had to assume that $a\in C(m),$ where
\begin{gather}
\label{Cm}
C(m)=\Big\{z\in\C; \; \Im(z)>0 \text{ and } 2\sqrt m\Im(z)\ge(1-m)|\Re(z)|\Big\}.
\end{gather}
Here and after, for $z\in\C,$ $\Re(z),$ $\Im(z)$ and $\ovl z$ denote the real part, the imaginary part and the conjugate of $z,$ respectively. Notice that in this paper, we have $m=0$ and $a\in C(0)$ in~\eqref{nlsg}. Hence $a=\vi\mu,$ for a positive real $\mu$ (the saturated case). Another set of complex numbers which plays an important role is:
\begin{gather}
\label{Dm}
D(m)=\Big\{z\in\C; \; \Im(z)>0 \text{ and } 2\sqrt m\Im(z)=(1-m)\Re(z)\Big\}.
\end{gather}
Note that $D(0)=C(0),$ $D(1)=\emptyset,$ and
\begin{align*}
& C(0)=\Big\{z\in\C; \; \Re(z)=0 \text{ and } \Im(z)>0\Big\},	\\
& C(1)=\Big\{z\in\C; \; \Im(z)>0\Big\}.
\end{align*}
Below, we summarize the results about the existence, uniqueness and finite time extinction property of the solutions in the previous literature according the presence of a potential $V(x)$ in the equation.

\begin{center}
Case without potential $(V\equiv0)$
\medskip

\renewcommand{\arraystretch}{1.2}
\begin{tabular}{|c|c|c|c|c|c|c|c|}
\hline
\multirow{3}{*}{$m$}	& \multirow{3}{*}{$a$}	&	\multirow{2}{*}{$\Omega\subseteq\R^N$}		&	\multicolumn{3}{c|}{Existence and uniqueness}	&
			\multirow{2}{*}{Finite time}		&	\multirow{3}{*}{Ref.}		\\
		&	&	\multirow{2}{*}{arbitrary}	&	\multicolumn{3}{c|}{of the solutions in}		&	\multirow{2}{*}{extinction for}	&		\\
\cline{4-6}
		&	&						&	$L^2(\Omega)$		&	$H^1_0(\Omega)$	&	$H^2(\Omega)$	&		&		\\
\hline
\multirow{2}{*}{$\in(0,1)$}	& 	\multirow{2}{*}{$\in C(m)$}		&	\multirow{2}{*}{$|\Omega|<\infty$}	&	\multirow{2}{*}{yes}
	&	\multirow{2}{*}{yes}	&	\multirow{2}{*}{yes}		&	$(N=1)$ or $(N\le3,$			&	\multirow{2}{*}{\cite{MR4053613}}	\\
	& 	&	&	&		&						&	$\Omega$ bounded and $C^{1,1})$	&							\\
\hline
$=1$		& $\in C(1)$	&	yes	&	yes	&	yes	&	yes	&	Impossible	&	\cite{MR4053613}						\\
\hline
$\in(0,1)$	& $\in C(m)\setminus D(m)$		&	$\Omega=\R^N$	&	yes	&	yes	&	yes	&	$N\le3$	&	\cite{MR4098330}	\\
\hline
\end{tabular}
\end{center}

\begin{center}
Case with a potential $V(x)$
\medskip

\renewcommand{\arraystretch}{1.2}
\begin{tabular}{|c|c|c|c|c|c|c|c|}
\hline
\multirow{3}{*}{$m$}	& \multirow{3}{*}{$a$}	&	\multirow{2}{*}{$\Omega\subseteq\R^N$}		&	\multicolumn{3}{c|}{Existence and uniqueness}	&
			\multirow{2}{*}{Finite time}		&	\multirow{3}{*}{Ref.}		\\
&	&	\multirow{2}{*}{arbitrary}	&	\multicolumn{3}{c|}{of the solutions in}		&	\multirow{2}{*}{extinction for}	&			\\
\cline{4-6}
					&				&			&	$L^2(\Omega)$		&	$H^1_0(\Omega)$	&	$H^2(\Omega)$		&
									&						\\
\hline
$\in(0,1)$	& $\in C(m)\setminus D(m)$		&	yes	&	yes	&	yes	&	yes	&	$N\le3$		&	\cite{MR4340780}	\\
\hline
$=1$		&	$\in C(1)$				& 	yes	&	yes	&	yes	&	yes	&	Impossible	&	\cite{MR4340780}	\\
\hline
$=0$		&	$\in C(0)$	&	$|\Omega|<\infty$	&	yes	&	yes	&	yes	&	$N\le3$		&	\cite{MR4340780}	\\
\hline
\multirow{2}{*}{$\in(0,1)$}	&   \multirow{2}{*}{$\in D(m)$}	&	\multirow{2}{*}{yes}	&	\multirow{2}{*}{yes}	&	\multirow{2}{*}{yes}	&
									yes		&	$(N=1)$ or $(N\le3$	&	\multirow{2}{*}{\cite{MR4503241}}	\\
	&		&		&		&		&	if $|\Omega|<\infty$		&	 and $|\Omega|<\infty)$	&			\\
\hline
$=0$		&	$\in C(0)$	&	yes			&	yes	&	yes		&	no	&	$N=1$			&	Here		\\
\hline
\end{tabular}
\end{center}

Due to the lack of regularity of the nonlinear term in equation \eqref{nls} (and in \eqref{nlsg} when $m\in[0,1)),$ a good technique to get the existence and uniqueness of solutions is to understand the equation as a special case of the abstract Cauchy problem
\begin{gather*}
\frac{\d u}{\d t}+Au=f,
\end{gather*}
where $Au=-\vi\Delta u-\vi Vu-\vi a|u|^{-(1-m)}u.$ In particular, we used the theory of maximal monotone operator in $L^2(\Omega)$ (see Brezis~\cite{MR0348562}).

The more favorable case corresponds to when $|\Omega|<\infty,$ $m\in(0,1]$ and $a\in C(m).$ So, in \cite{MR4053613} we show directly that $(D(A),A)$ is maximal monotone by using the embedding $L^p(\Omega)\inj L^2(\Omega),$ for any $p\ge2:$ $\left||u|^{-(1-m)}u\right|=|u|^m\in L^\frac2m(\Omega)\inj L^2(\Omega),$ $\forall m\in(0,1].$

When $\Omega=\R^N$ but $m\in(0,1)$ and $a\in C(m)\setminus D(m)$ the abstract theory can be applied. So, in \cite{MR4098330} we show that $(D(A),A)$ is maximal monotone in the following way. First, we build solutions compactly supported in $H^2(\R^N)$ to $(A+I)u=F$ with help of the results in \cite{MR2876246,MR3315701}. Second, we obtain a priori estimates in $H^2(\R^N)$ with \cite[Lemma~4.2]{MR4098330}. Third, we show that $(D(A),A)$ is maximal monotone by approximations with solutions compactly supported.

An approximation argument allowed us to extend the above mentioned results to the case in which $m\in[0,1],$ $a\in C(m)\setminus D(m)$ but according the boundedness or not of $\Omega$ (\cite{MR4340780}). First, we approximate $(D(A),A)$ by a nice maximal monotone operator $(D(A),A_\eps).$ Second, we obtain a priori estimates in $H^2(\Omega)$ with \cite[Lemma~4.2]{MR4098330}. Third, we pass to the limit. If $m=0$ then it is assumed that $|\Omega|<\infty,$ and then $\left|\frac{u}{|u|}\right|=1\in L^\infty(\Omega)\inj L^2(\Omega).$

The critical case $a\in D(m)$ was considered in \cite{MR4503241} when $m\in(0,1).$ The main idea was to get different a priori estimates for the approximate problems, with $(D(A),A_\eps).$

One of the main difficulties when considering the problem \eqref{nls}, main goal of this article, comes from the lack of separability of the space $L^\infty(\Omega)$ (notice that the nonlinear term is a bounded function). When $|\Omega|<\infty$ it is possible to use additional a priori estimates allowing the passing to the limit when $\eps\searrow0,$ but this strategy fails when $\Omega$ is unbounded.

We point out that although there are some abstract results (and with many interesting applications to nonlinear PDEs -- see Bénilan and Ha~\cite{MR0405187}) this theory does not apply to the case of unbounded domains $\Omega.$ The difficulty comes from the fact that it is very hard to get ``strong solutions" (i.e., solutions such that $u_t\in L^1(0,T;H^{-1}(\Omega)+L^\infty(\Omega)),$ for any $T>0).$ This is also associated to the fact that when the Banach space $X$ is not reflexive (as it is also the case of $X=H^1_0(\Omega)\cap L^1(\Omega))$ then the Radon-Nikod\'ym property may fail (see, e.g., Diestel and Uhl~\cite{MR0453964}). Actually, we do not even know if the time differential $u_t$ is measurable when the data $(u_0$ and $f)$ are not regular enough. Nevertheless, by introducing a suitable sequence of spaces $(Y_n)_{n\in\N}$ satisfying (in an appropriate sense) the Radon-Nikod\'ym property (see Definition~\ref{defRNP} below) it is possible to get a limit solution which is stronger than the ``weak solutions" but weaker than the ``strong solutions".

The paper is organized as follows. Section~\ref{exiuni} is devoted to state the existence and uniqueness of weak and strong solutions under suitable assumptions on the initial data $u_0$ and $f.$ The explanation at the sense in which the weak solutions satisfy \eqref{nls} already uses the approximating sequence of spaces $(Y_n)_{n\in\N}$ mentioned before. Section~\ref{finite} contains the statement of the strong stabilization results when $t\nearrow\infty.$ In particular, when $N=1$ we get the finite time extinction of solutions with a time-decay estimate that is stronger than the usual time-decay estimate for sublinear parabolic equations (see, e.g., \cite{MR2002i:35001}). The occurrence of the problem of measurability of $u_t$ is presented in Section~\ref{secnm} (for a different approach in an abstract framework, see Deville~\cite{MR1045133}). The proof of the existence and uniqueness results stated in Section~\ref{exiuni} are presented in Section~\ref{proofexi}, after recalling and developing some results on Functional Analysis in 
Section~\ref{secfunana}. The results stated in Section~\ref{finite} are proved in Section~\ref{proofext}. The consideration of solutions with lower regularity is carried out in Section~\ref{secaH1} and explains how it is possible to pass to the limit in suitable approximations.

To end this introduction, we collect here some notations which will be used throughout this paper. For a real number $t\in\R,$ $t_+=\max\{t,0\}$ is its positive part. Let $\Omega$ be an open subset of $\R^N.$ Unless specified, all functions are complex-valued $(H^1(\Omega)\eqdef H^1(\Omega;\C),$ etc) and all the vector spaces are considered over the field $\R.$ For $p\in[1,\infty],$ $p^\prime$ is the conjugate of $p$ defined by $\frac{1}{p}+\frac{1}{p^\prime}=1.$ For a (real) Banach space $X,$ we denote by $X^\star\eqdef\Lr(X;\R)$ its topological dual and by $\langle\: . \; , \: . \:\rangle_{X^\star,X}\in\R$ the $X^\star-X$ duality product. When $X$ (respectively, $X^\star)$ is endowed of the weak topology $\sigma\left(X,X^\star\right)$ \big(respectively, the weak$\star$ topology $\sigma\left(X^\star,X\right)\big),$ it is denoted by $X_\w$ (respectively, by $X_{\w\star}).$ For $p\in[1,\infty],$ $u\in L^p_\loc\big([0,\infty);X\big)$ means that $u\in L^p_\loc(0,\infty;X)$ and for any $T>0,$ $u_{|(0,T)}\in L^p(0,T;X).$ In the same way, we will use the notation $u\in W^{1,p}_\loc\big([0,\infty);X\big).$ The scalar product in $L^2(\Omega)$ between two functions $u,v$ is, $(u,v)_{L^2(\Omega)}=\Re\int_{\Omega}u(x)\ovl{v(x)}\d x.$ For any real interval $I$ and Banach space $X,$ $C_\w(I;X)$ is the space of weakly continuous function from $I$ to $X_\w.$ The space of measurable functions $u:\Omega\tends\C$ such that $|u|<\infty,$ almost everywhere in $\Omega,$ is denoted by $L^0(\Omega).$ Auxiliary positive constants will be denoted by $C$ and may change from a line to another one. Also for positive parameters $a_1,\ldots,a_n,$ we shall write $C(a_1,\ldots,a_n)$ to indicate that the constant $C$ depends only and continuously on $a_1,\ldots,a_n.$ The set of positive integers is denoted by $\N,$ and $\N_0=\N\cup\{0\}.$ We denote by $\ovl B_{L^\infty}(0,1)$ the closed unit ball of $L^\infty(\Omega).$ Finally, the Lebesgue measure of a measurable set $A\subset\R^N$ will be denoted by $|A|.$

\section{Existence and uniqueness of weak and strong solutions}
\label{exiuni}

Throughout this paper we shall always identify $L^2(\Omega)$ with its topological dual. We refer to Section~\ref{secfunana} for some results of Functional Analysis. The following assumptions will be needed to build solutions.
\begin{assum}
\label{ass1}
We assume that
\begin{gather}
\label{O}
\Omega \text{ is any nonempty open subset of } \R^N,	\\
\label{a}
\mu>0,										\\
\label{V}
V\in L^\infty(\Omega;\R)+L^{p_V}(\Omega;\R),
\end{gather}
where,
\begin{gather}
\label{pV}
p_V=
\begin{cases}
2,							&	\text{if } N=1,	\\
2+\beta, \text{ for some } \beta>0,	&	\text{if } N=2,	\\
N,							&	\text{if } N\ge3.
\end{cases}
\end{gather}
\medskip
\end{assum}

\noindent
Now, let us recall the definition of solution (\cite{MR4340780}) for a general $\mu\in\C.$

\begin{defi}
\label{defsol}
Assume~\eqref{O}, \eqref{V} and \eqref{pV}. Let $\mu\in\C,$ $f\in L^1_\loc\big([0,\infty);L^2(\Omega)\big)$ and $u_0\in L^2(\Omega).$ We shall say that $u$ is a \textit{strong solution} (or an $H^1_0$-\textit{solution}) to \eqref{nls}--\eqref{u0} if $u$ satisfies the following properties.
\begin{enumerate}
\item
\label{defsol1}
We have that,
\begin{gather}
\label{defsol11}
u\in L^1_\loc\big([0,\infty);X\big)\cap W^{1,\infty}_\loc\big([0,\infty);X^\star\big)\inj C\big([0,\infty);L^2(\Omega)\big),
\end{gather}
where $X=H^1_0(\Omega)\cap L^1(\Omega)$ (hence, $X^\star=H^{-1}(\Omega)+L^\infty(\Omega)).$
\item
\label{defsol2}
There exists a \textit{saturated section} $U$ associated to $u,$ namely a $U\in L^\infty\big((0,\infty)\times\Omega\big)$ such that
\begin{gather}
\label{defsol21}
\vi\frac{\partial u}{\partial t}+\Delta u+V(x)u+\vi\mu\,U=f(t,x), \text{ in } \Dr^\p\big((0,\infty)\times\Omega\big),	\\
\label{defsol22}
\|U\|_{L^\infty((0,\infty)\times\Omega)}\le1,														\\
\label{defsol23}
U(t,x)=\frac{u(t,x)}{|u(t,x)|}, \text{ as soon as } u(t,x)\neq0.
\end{gather}
\item
\label{defsol3}
$u(0)=u_0,$ in $L^2(\Omega).$
\end{enumerate}
We shall say that $u$ is a \textit{weak solution} (or an $L^2$-\textit{solution}) to \eqref{nls}--\eqref{u0} if there exists,
\begin{gather}
\label{fn}
(u_n,U_n,f_n)_{n\in\N}\subset C\big([0,\infty);L^2(\Omega)\big)\times L^\infty\big((0,\infty)\times\Omega\big)\times L^1_\loc\big([0,\infty);L^2(\Omega)\big),
\end{gather}
such that for any $n\in\N,$ $u_n$ is a strong solution of \eqref{nls}--\eqref{nlsb} with the saturated section $U_n,$ where the right-hand side member of \eqref{nls} is $f_n,$ and if
\begin{gather}
\label{cv}
(u_n,f_n)\xrightarrow[n\tends\infty]{C([0,T];L^2(\Omega))\times L^1(0,T;L^2(\Omega))}(u,f),	\\
\label{U}
U_n\underset{n\to\infty}{\overset{L^\infty((0,T)\times\Omega)_{\w\star}}{-\!\!\!-\!\!\!-\!\!\!-\!\!\!-\!\!\!-\!\!\!-\!\!\!-\!\!\!-\!\!\!-\!\!\!-\!\!\!\weak}}U,
\end{gather}
for any $T>0.$ Sometimes, we shall write $(u,f),$ $(u,U),$ or $(u,U,f)$ to designate a solution with the obvious meanings.
\end{defi}

\begin{rmk}
\label{rmkdefsol}
The embedding \eqref{defsol11} comes from Theorem~\ref{thmembL2} below.
\medskip
\end{rmk}

\noindent
Before recalling a result of uniqueness and continuous dependence (\cite[Proposition~2.5]{MR4340780}), we explain below in which way the weak solutions satisfy~\eqref{nls}.

\begin{prop}
\label{propsolL2}
Assume~\eqref{O}, \eqref{V} and \eqref{pV}. Let $\mu\in\C,$ $f\in L^1_\loc\big([0,\infty);L^2(\Omega)\big)$ and $u_0\in L^2(\Omega).$ Let $(Y_n)_{n\in\N_0}$ be any $L^1$-approximating sequence of RNP-spaces $($see Definition~$\ref{defRNP}$ below$).$ If $u$ is a weak solution to \eqref{nls}--\eqref{u0} then for any $n\in\N_0,$
\begin{gather}
\label{propsolL21}
u\in W^{1,1}_\loc\big([0,\infty);H^{-2}(\Omega)+Y^\star_n\big),
\end{gather}
and $u$ solves~\eqref{nls} in $L^1_\loc\big([0,\infty);H^{-2}(\Omega)+Y^\star_n\big)\inj\Dr^\p\big((0,\infty)\times\Omega\big).$ Finally there exists $N_0\subset(0,\infty)$ with $|N_0|=0$ such that,
\begin{gather}
\label{propsolL22}
u^\p(t)\in H^{-2}(\Omega)+L^\infty(\Omega),
\end{gather}
for any $t\in(0,\infty)\cap N_0^\co.$ In particular, $u$ solves~\eqref{nls} in $ H^{-2}(\Omega)+L^\infty(\Omega),$ for almost every $t>0.$
\end{prop}

\begin{rmk}
\label{rmkpropsolL2}
Below are some comments about Proposition~\ref{propsolL2}.
\begin{enumerate}
\item
\label{rmkpropsolL21}
If $|\Omega|<\infty$ then the spaces $H^{-2}(\Omega)+Y^\star_n$ may be replaced with $H^{-2}(\Omega).$ See the end of the proof of Proposition~\ref{propsolL2} for the details.
\item
\label{rmkpropsolL22}
``RNP" stands for ``Radon-Nikod\'ym Property". For the justification of this terminology, see Section~\ref{secnm} below.
\item
\label{rmkpropsolL23}
Whether or not $u^\p:[0,\infty)\tends H^{-2}(\Omega)+L^\infty(\Omega)$ is measurable is an open question.
\end{enumerate}
\end{rmk}

\begin{prop}[\textbf{Uniqueness and continuous dependence}]
\label{propdep}
Let Assumption~$\ref{ass1}$ be fulfilled, let $f,\wt f\in L^1_\loc\big([0,\infty);L^2(\Omega)\big)$ and $X=H^1_0(\Omega)\cap L^1(\Omega).$ Finally, let $p\in[1,\infty]$ and let
\begin{gather}
\label{propdephyp}
u,\wt u\in L^p_\loc\big([0,\infty);X\big)\cap W^{1,p^\p}_\loc\big([0,\infty);X^\star\big)\inj C\big([0,\infty);L^2(\Omega)\big),
\end{gather}
be solutions in $\Dr^\p\big((0,\infty)\times\Omega\big)$ to,
\begin{gather*}
\vi u_t+\Delta u+Vu+\vi\mu\,U=f,				\\
\vi\wt {u_t}+\Delta\wt u+V\wt u+\vi\mu\,\wt U=\wt f ,
\end{gather*}
respectively, where $U$ and $\wt U$ satisfy~\eqref{defsol22}--\eqref{defsol23}. Then,
\begin{gather}
\label{estthmweak}
\|u(t)-\wt u(t)\|_{L^2(\Omega)}\le\|u(s)-\wt u(s)\|_{L^2(\Omega)}+\vint_s^t\|f(\sigma)-\wt f(\sigma)\|_{L^2(\Omega)}\d\sigma,
\end{gather}
for any $t\ge s\ge0.$ Finally,~\eqref{estthmweak} also holds true for the weak solutions.
\end{prop}

\begin{thm}[\textbf{Existence and uniqueness of weak solutions}]
\label{thmweak}
Let Assumption~$\ref{ass1}$ be fulfilled and let $f\in L^1_\loc\big([0,\infty);L^2(\Omega)\big).$ Then for any $u_0\in L^2(\Omega),$ there exists a unique weak solution $u$ to \eqref{nls}--\eqref{u0}. In addition,
\begin{gather}
\label{L1}
u\in L^1_\loc\big([0,\infty);L^1(\Omega)\big),	\\
\label{L2+}
\frac12\|u(t)\|_{L^2(\Omega)}^2+\mu\dsp\vint_s^t\|u(\sigma)\|_{L^1(\Omega)}\d\sigma
			\le\frac12\|u(s)\|_{L^2(\Omega)}^2+\Im\dsp\iint\limits_{s\;\Omega}^{\text{}\;\;t}f(\sigma,x)\,\ovl{u(\sigma,x)}\,\d x\,\d\sigma,
\end{gather}
for any $t\ge s\ge0.$ If $|\Omega|<\infty$ then the inequality in~\eqref{L2+} becomes an equality.
\end{thm}

\begin{rmk}
\label{rmkthmweak}
If $(u,f)$ and $(\wt u,\wt f)$ are two weak solutions then, by Hölder's inequality, we obtain for any $p\in(1,2),$
\begin{gather*}
\|u-\wt u\|_{L^\frac{p}{2-p}(s,t;L^p(\Omega))}\le\|u-\wt u\|_{L^1(s,t;L^1(\Omega))}^\frac{2-p}p\|u-\wt u\|_{C_\b([s,t];L^2(\Omega))}^\frac{2(p-1)}p,
\end{gather*}
for any $t\ge s\ge0.$ Then by \eqref{estthmweak}--\eqref{L2+}, the left-hand side of the above estimate is bounded.
\end{rmk}

\begin{thm}[\textbf{Existence and uniqueness of strong solutions}]
\label{thmstrongaH2}
Let Assumption~$\ref{ass1}$ be fulfilled and $f\in W^{1,1}_\loc\big([0,\infty);L^2(\Omega)\big).$ Then for any $u_0\in H^2_\loc(\Omega)\cap H^1_0(\Omega)\cap L^1(\Omega)$ for which
\begin{gather*}
\left(\Delta u_0+\vi\mu\frac{u_0}{|u_0|}\right)_{|\{u_0\neq0\}}\in L^2(\{u_0\neq0\}),
\end{gather*}
there exists a unique $H^1_0$-solution $u$ to \eqref{nls}--\eqref{u0}. Furthermore, $u$ satisfies~\eqref{nls} in $L^\infty_\loc\big([0,\infty);L^2_\loc(\Omega)\big)\cap L^\infty_\loc\big([0,\infty);H^{-1}(\Omega)\big)$ as well as the following properties.
\begin{enumerate}
\item
\label{thmstrongaH21}
$u\in C_\w\big([0,\infty);H^1_0(\Omega)\big)\cap L^\infty_\loc\big([0,\infty);L^1(\Omega)\big)\cap W^{1,\infty}_\loc\big([0,\infty);L^2(\Omega)\big).$
\item
\label{thmstrongaH22}
For any $t\ge s\ge0,$
\begin{empheq}[left=\empheqlbrace]{align}
\label{strongaH21}
	&	\|u(t)-u(s)\|_{L^2(\Omega)}\le\|u_t\|_{L^\infty(s,t;L^2(\Omega))}|t-s|,		\frac{}{} \\
\label{strongaH22}
	&	\|u(t)\|_{L^2(\Omega)}\le A(t),									\frac{}{} \\
\label{strongaH23}
	&	\left\|u_t\right\|_{L^\infty(0,t;L^2(\Omega))}\le B(t),					\frac{}{} \\
\label{strongaH24}
	&	\|\nabla u(t)\|_{L^2(\Omega)}^2+\mu\|u(t)\|_{L^1(\Omega)}\le C(t)A(t),
\end{empheq}
where,
\begin{align*}
&	A(t)=\|u_0\|_{L^2(\Omega)}+\int_0^t\|f(s)\|_{L^2(\Omega)}\d s,										\\
&	B(t)=\|\Delta u_0+Vu_0+\vi\mu U_0-f(0)\|_{L^2(\Omega)}+\int_0^t\|f^\p(\sigma)\|_{L^2(\Omega)}\d\sigma,		\\
&	C(t)=C\left(A(t),B(t),\|f(t)\|_{L^2(\Omega)},\|V_1\|_{L^\infty(\Omega)},\|V_2\|_{L^{p_V}(\Omega)},N,\beta\right),
\end{align*}
for some $U_0\in \ovl B_{L^\infty}(0,1)$ satisfying $U_0=\frac{u_0}{|u_0|},$ almost everywhere where $u_0\neq0.$
\item
\label{thmstrongaH23}
The map $t\longmapsto\|u(t)\|_{L^2(\Omega)}^2$ belongs to $W^{1,\infty}_\loc\big([0,\infty);\R\big)$  and we have,
\begin{gather}
\label{L2}
\frac12\frac{\d}{\d t}\|u(t)\|_{L^2(\Omega)}^2+\mu\|u(t)\|_{L^1(\Omega)}=\Im\vint_{\Omega}f(t,x)\,\ovl{u(t,x)}\,\d x,
\end{gather}
for almost every $t>0.$
\item
\label{thmstrongaH24}
If $f\in W^{1,1}(0,\infty;L^2(\Omega))$ then $u\in L^\infty(0,\infty;H^1_0(\Omega)\cap L^1(\Omega))\cap W^{1,\infty}(0,\infty;L^2(\Omega)).$
\end{enumerate}
\end{thm}

\begin{rmk}
\label{rmkf0}
Below are some comments about Theorem~\ref{thmstrongaH2}.
\begin{enumerate}
\item
\label{rmkf01}
The solution obtained in Theorem~\ref{thmstrongaH2} is not an $H^2$-solution (in the sense of~\cite[Definition 2.2]{MR4340780}). Indeed, to be one, we would have to have $\Delta u(t)\in L^2(\Omega),$ for almost every $t>0,$ while we merely have $\Delta u(t)\in L^2_\loc(\Omega),$ for almost every $t>0.$ 
Existence, uniqueness, finite time extinction and asymptotic behavior of the $H^2$-solutions are obtained in \cite{MR4340780} in the special case of $|\Omega|<\infty.$
\item
\label{rmkf03}
Using that $u\in C^{0,1}\big([0,\infty);L^2(\Omega)\big),$ and the Gagliardo-Nirenberg and Hölder inequalities, we get that for any $p\in\left(1,\frac{2N}{N-2}\right)$ $(p\in(1,\infty]$ if $N=1),$
\begin{gather*}
u\in C^{0,\alpha}\big([0,\infty);L^p(\Omega)\big)	\;
\big(u\in C^{0,\alpha}_\b\big([0,\infty);L^p(\Omega)\big), \text{ if } f\in W^{1,1}(0,\infty;L^2(\Omega))\big),
\end{gather*}
where $\alpha=\frac{2N-p(N-2)}{2p}$ if $p\ge2,$ and $\alpha=2\frac{p-1}p$ if $p\le2.$
\end{enumerate}
\end{rmk}

\section{Finite time extinction and asymptotic behavior for large time}
\label{finite}

\section*{Large time behavior of the weak solutions}

\begin{thm}
\label{thm0w}
Let Assumption~$\ref{ass1}$ be fulfilled, $f\in L^1(0,\infty;L^2(\Omega)),$ $u_0\in L^2(\Omega)$ and let $u$ be the unique weak solution to \eqref{nls}--\eqref{u0} given by Theorem~$\ref{thmweak}.$ Then,
\begin{gather*}
\vlim_{t\nearrow\infty}\|u(t)\|_{L^2(\Omega)}=0.
\end{gather*}
\end{thm}

\section*{Finite time extinction and asymptotic behavior of the strong solutions}

Under some additional conditions on $f(t,x)$ we have:
\begin{thm}[\textbf{Finite time extinction and time decay estimates}]
\label{thmextN1}
Let Assumption~$\ref{ass1}$ be fulfilled, $f\in W^{1,1}(0,\infty;L^2(\Omega)),$ $u_0\in H^2_\loc(\Omega)\cap H^1_0(\Omega)\cap L^1(\Omega)$ with
\begin{gather*}
\left(\Delta u_0+\vi\mu\frac{u_0}{|u_0|}\right)_{|\{u_0\neq0\}}\in L^2(\{u_0\neq0\}),
\end{gather*}
and let $u$ be the unique strong solution to \eqref{nls}--\eqref{u0} given by Theorem~$\ref{thmstrongaH2}.$
\begin{enumerate}
\item
\label{thmextT*}
If there exists $T_0\in[0,\infty)$ such that
\begin{gather}
\label{f}
f\in L^\infty\big((T_0,\infty)\times\Omega\big) \, \text{ and } \; \|f\|_{L^\infty((T_0,\infty)\times\Omega)}<\mu,
\end{gather}
then the following holds true.
\begin{itemize}
\item
\label{thmextN11}
If $N=1$ then
\begin{gather}
\label{01}
\forall t\ge T_\star, \; \|u(t)\|_{L^2(\Omega)}=0,
\end{gather}
for some,
\begin{gather}
\label{T*N1}
T_0\le T_\star\le C\|u(T_0)\|_{L^2(\Omega)}^\frac12\|\nabla u\|_{L^\infty(0,\infty;L^2(\Omega))}^\frac12+T_0,
\end{gather}
where $C=C(\mu-\|f\|_{L^\infty((T_0,\infty)\times\Omega)}).$ Actually,
\begin{gather*}
\|u(t)\|_{L^2(\Omega)}\le\left(\|u_0\|_{L^2(\Omega)}^\frac12+C(T_0-t)\right)_+^2,
\end{gather*}
for any $t\ge T_0$ and for some $C=C(\|\nabla u\|_{L^\infty(0,\infty;L^2(\Omega))},\mu-\|f\|_{L^\infty((T_0,\infty)\times\Omega)}).$
\item
\label{thmextN12}
If $N=2$ then for any $t\ge T_0,$
\begin{gather}
\label{thmrtdH11}
\|u(t)\|_{L^2(\Omega)}\le\|u(T_0)\|_{L^2(\Omega)}e^{-C(t-T_0)},
\end{gather}
where $C=C(\|\nabla u\|_{L^\infty(0,\infty;L^2(\Omega))},\mu-\|f\|_{L^\infty((T_0,\infty)\times\Omega)}).$
\item
\label{thmextN13}
If $N\ge3$ then for any $t\ge T_0,$
\begin{gather}
\label{thmrtdH12}
\|u(t)\|_{L^2(\Omega)}\le\dfrac{\|u(T_0)\|_{L^2(\Omega)}}
{\left(1+C\|u(T_0)\|_{L^2(\Omega)}^\frac{(N-2)}2(t-T_0)\right)^\frac2{(N-2)}},
\end{gather}
where $C=C(\|\nabla u\|_{L^\infty(0,\infty;L^2(\Omega))},\mu-\|f\|_{L^\infty((T_0,\infty)\times\Omega)},N).$
\end{itemize}
\item
\label{thmextN14}
Suppose $N=1.$ For each $M>0,$ there exists $\eps_\star=\eps_\star(M,\mu)$ satisfying the following property. Let $T_0\ge0.$ If
\begin{gather}
\label{thmextN1e1}
\begin{cases}
\|f\|_{W^{1,1}(0,\infty;L^2(\Omega))}\le M,													\medskip \\
\|\nabla u_0\|_{L^2(\Omega)}+\left\|\Delta u_0+\vi\mu\frac{u_0}{|u_0|}\right\|_{L^2(\{u_0\neq0\})}\le M,	\medskip \\
\|u_0\|_{L^2(\Omega)}+\|f\|_{L^1(0,\infty;L^2(\Omega))}\le\eps_\star,								\medskip \\
\|u_0\|_{L^2(\Omega)}\le\eps_\star T_0^2,													\medskip \\
\|f(t)\|_{L^2(\Omega)}\le\eps_\star\big(T_0-t\big)_+, \; \forall t\ge0,
\end{cases}
\end{gather}
then~\eqref{01} holds true with $T_\star=T_0.$
\end{enumerate}
\end{thm}

\begin{rmk}
\label{rmkthmextN1}
If $f$ satisfies~\eqref{f} then it follows from~\eqref{01} and~\eqref{nls} that,
\begin{gather*}
\vi\mu\,U(t,x)=f(t,x),
\end{gather*}
for almost every $(t,x)\in(T_\star,\infty)\times\Omega.$
\end{rmk}

\section{On the non-measurability}
\label{secnm}

Due to the special nonlinearity, in this paper, we are not able to build $H^1_0$-solutions of~\eqref{nls} under the mere ``natural conditions"
\begin{gather}
\label{secnmu0}
u_0\in H^1_0(\Omega),	\\
\label{secnmf}
f\in L^1_\loc\big([0,\infty);H^1_0(\Omega)\big)\cap L^\infty_\loc\big([0,\infty);H^{-1}(\Omega)+L^\infty(\Omega)\big),
\end{gather}
(see Theorem~\ref{thmstrongaH2}). In the previous papers~\cite{MR4098330,MR4053613,MR4340780,MR4503241}, the methods used to obtain such solutions have been the following. We proceed by density, either by starting from the equation~\eqref{nlsg}, or by starting from the approximate equation
\begin{gather}
\label{nlsgm}
\vi\frac{\partial u}{\partial t}+\Delta u+V(x)u+a\frac{u}{(|u|^2+\eps)^\frac{1-m}2}=f(t,x), \text{ in } (0,\infty)\times\Omega,
\end{gather}
where the nonlinearity $|u|^{-(1-m)}u$ is regularized as above. In both cases, we show that the sequence of solutions $(u_\eps)_{\eps>0}$ is bounded in $L^{m+1}(I;H^1_0(\Omega)\cap L^{m+1}(\Omega))$ and in
\begin{gather}
\label{Wm}
W^{1,\frac{m+1}m}(I;H^{-1}(\Omega)+L^\frac{m+1}m(\Omega)),
\end{gather}
for any $T>0,$ where $I=(0,T).$ When $0<m\le1$ then $1<\frac{m+1}m<\infty$ and thus $H^{-1}(\Omega)+L^\frac{m+1}m(\Omega)$ is reflexive, and so is the space of \eqref{Wm}. It follows that we may extract from $(u_\eps)_{\eps>0}$ a sequence which converges in the weak topology to some $u$ belonging to \eqref{Wm}. When $m=0,$ in \cite{MR4340780} we assumed that $|\Omega|<\infty.$ It follows that $L^\infty(\Omega)\inj H^{-1}(\Omega),$ and the space in \eqref{Wm} becomes 
\begin{gather}
\label{W0b}
W^{1,\infty}(I;H^{-1}(\Omega)).
\end{gather}
But since $H^1_0(\Omega)$ is reflexive and separable then $L^1(I;H^1_0(\Omega))$ is separable and is the predual of $L^\infty(I;H^{-1}(\Omega))$ (see Section~\ref{secfunana} below). And again we may extract from $(u_\eps)_{\eps>0}$ a sequence which converges in the weak$\star$ topology to some $u$ belonging to the space in \eqref{W0b}. But in the present paper, $m=0$ and it is no more assumed that $|\Omega|<\infty.$ Then, we have to deal with the space
\begin{gather}
\label{W0}
W^{1,\infty}(I;X^\star),
\end{gather}
where $X=H^1_0(\Omega)\cap L^1(\Omega).$ But $X$ is not reflexive so that $L^1(I;X)$ is not the predual of $L^\infty(I;X^\star).$ We then have to proceed in a different way. Using the Arzelà-Ascoli compactness Theorem, we could extract a sequence $(u_{\eps_n})_{n\in\N}$ of $(u_\eps)_{\eps>0}$ such that for any $t\in\ovl I,$ $(u_{\eps_n}(t))_{n\in\N}$ converges to some $u\in L^\infty(I;H^1_0(\Omega))$ in the weak topology $H^1_0(\Omega)_\w.$ But, $H^1_0(\Omega)\inj X^\star,$ whose the predual $X$ is separable. We may deduce that,
\begin{gather}
\label{ulip}
\|u(t)-u(s)\|_{X^\star}\le\liminf_{n\to\infty}\|u_{\eps_n}(t)-u_{\eps_n}(s)\|_{X^\star}\le C|t-s|,
\end{gather}
for any $t,s\in\ovl I.$ In particular, $u:I\tends X^\star$ is absolutely continuous. But it is well-known that $X^\star$ does not satisfy the Radon-Nikod\'ym property, RNP in short (Stegall~\cite{MR374381}, Maurey~\cite{MR0402503}), and even not the weak Radon-Nikod\'ym property, WRNP in short (Musia\l~\cite{MR537118}, Rosenthal~\cite{MR358307}), so that we cannot conclude that $u:I\tends X^\star$ is differentiable almost everywhere (Diestel and Uhl~\cite[Theorem~2, p.107]{MR0453964}). So, the idea is to embed $X^\star$ in a bigger space, but the smallest bigger space as possible, satisfying the RNP. For instance, this space is chosen to be reflexive (Phillips~\cite{MR4094}, Diestel and Uhl~\cite[Corollary~13, p.76]{MR0453964}), or to be a separable dual space (Dunford and Pettis~\cite[Theorem~2.1.4, p.345]{MR2020}, Diestel and Uhl~\cite[Corollary~1, p.79]{MR0453964}). In this case, for such a space $Y^\star,$ we obtain that $u\in C^{0,1}(\ovl I;Y^\star)=W^{1,\infty}(I;Y^\star).$ This justifies the Definition~\ref{defRNP} below. Nevertheless, we may wonder if, thereafter, we may give up this space. Indeed, from \eqref{nls} and the fact that $u\in W^{1,\infty}(I;Y^\star),$ there exist $N_0\subset I$ and $C>0$ such that $|N_0|=0,$ and
\begin{gather}
\label{u'}
\forall t\in I\setminus N_0, \; u^\p(t)\in X^\star \text{ and } \|u^\p(t)\|_{X^\star}\le C,			\\
\label{diffY*}
\forall t\in I\setminus N_0, \; \lim_{h\to0}\left\|\frac{u(t+h)-u(t)}h-u^\p(t)\right\|_{Y^\star}=0,	\\
\label{int}
\forall t,s\in\ovl I, \, u(t)-u(s)\overset{Y^\star}{=}\int_s^tu^\p(\sigma)\d\sigma.
\end{gather}
But $u\in C^{0,1}(\ovl I;X^\star),$ so that the equality in \eqref{int} also makes sense in $X^\star.$ From \eqref{ulip} and \eqref{diffY*}, if $Y\inj X$ with dense embedding (which will be our case), then we easily deduce that,
\begin{gather}
\label{ufm}
\forall t\in I\setminus N_0, \; \frac{u(t+h)-u(t)}h\underset{h\to0}{-\!\!\!-\!\!\!-\!\!\!-\!\!\!\weak}u^\p(t), \text{ in the weak}\!\star \text{topology } \sigma(X^\star,X).
\end{gather}
It follows that $u^\p:I\tends X^\star$ is weak$\star$-measurable, and by Pettis' Theorem (\cite{MR1501970}), it will be measurable if, and only if, $u^\p(I\setminus N_1)$ is separable in $X^\star,$ for some null set $N_1\subset I.$ The answer is no, as shows the counter-example below. Indeed, the separability of a subset of $L^\infty(\Omega)$ is very difficult to obtain in a general way. So, we build a solution $u$ with lower regularity, in the sense that we do not know whether $u^\p:I\tends H^{-1}(\Omega)+L^\infty(\Omega)$ is measurable (see Section~\ref{secaH1} below). For a related analysis, see Deville~\cite{MR1045133}.

\begin{exa}
\label{exanm}
Below, we give an example of a function $u\in L^\infty(\R;L^\infty(\R))$ and of a separable Hilbert space $\Sigma,$ such that $\Sigma\inj L^1(\R^N)$ with dense embedding, and which satisfy,
\begin{gather}
\label{exanm1}
u:\R\tends L^\infty(\R) \text{ is weakly}\!\star \, \text{differentiable everywhere,}	\\
\label{exanm2}
u:\R\tends\Sigma^\star \text{ is differentiable almost everywhere,}			\\
\label{exanm3}
\forall t\in\R, \; u^\p(t)\in L^\infty(\R) \text{ and } \|u^\p(t)\|_{L^\infty(\R)}\le1,		\\
\label{exanm4}
\forall t,s\in\R, \, u(t)-u(s)\overset{L^\infty(\R)}{=}\int_s^tu^\p(\sigma)\d\sigma,
\end{gather}
but
\begin{gather}
\label{exanm5}
u^\p:\R\tends L^\infty(\R) \text{ is not measurable.}
\end{gather}
We define $u:\R\tends L^\infty(\R)$ as follows. Let $t\in\R.$ For any $x\in\R,$ let $u(t)(x)=\arctan|t+x|.$ Clearly, $u\in C^{0,1}_\b(\R;C_\b(\R)).$ Let $t\in\R.$ Let $\vphi\in L^1(\R).$ By the dominated convergence Theorem, we have
\begin{gather*}
\lim_{h\to0}\left\langle\frac{u(t+h)-u(t)}h,\vphi\right\rangle_{L^\infty(\R),L^1(\R)}=\int_\R\frac{\sign(t+x)}{1+(t+x)^2}\vphi(x)\d x.
\end{gather*}
Therefore, \eqref{exanm1} and \eqref{exanm3} hold true. It follows from Pettis' Theorem (Diestel and Uhl~\cite[Corollary~4, p.42-43]{MR0453964}) that $u^\p:\R\tends L^\infty(\R)$ is measurable if, and only if, there exists a null set $N\subset\R$ such that $u^\p(\R\setminus N)$ is separable in $L^\infty(\R).$ Let $t,s\in\R$ with $t>s.$ If $t-s\ge2$ then we let $x_0=1-t,$ and we obtain $u^\p(t)(x_0)-u^\p(s)(x_0)\ge\frac12,$ while if $t-s<2$ then we let $x_0=-\frac{t+s}2,$ and we obtain $u^\p(t)(x_0)-u^\p(s)(x_0)\ge1.$ It follows that for any $t,s\in\R$ with $t\neq s,$ we have
\begin{gather*}
\|u^\p(t)-u^\p(s)\|_{L^\infty(\R)}\ge\frac12.
\end{gather*}
As a consequence, for any null set $N\subset\R,$ $u^\p(\R\setminus N)$ is not separable. Hence, \eqref{exanm5}. In particular, $u^\p\not\in L^\infty(\R;L^\infty(\R)).$ Nevertheless, if we consider,
\begin{gather*}
\Sigma=\Big\{u\in L^2(\R); |\:.\:|u(\:.\:)\in L^2(\R)\Big\},
\end{gather*}
with its obvious norm, then it is well-known that $\Sigma$ is a separable Hilbert space and that $\Sigma\inj L^1(\R)$ with dense embedding. Therefore, $\Sigma^\star$ is reflexive and $L^\infty(\R)\inj\Sigma^\star$ with dense embedding. Then $\Sigma^\star$ has the Radon-Nikod\'ym property, and since $u:\R\tends\Sigma^\star$ is Lipschitz continuous, \eqref{exanm2} follows. In particular, $u^\p\in L^\infty(\R;\Sigma^\star)$ and $u^\p:\R\tends\Sigma^\star$ is Bochner integrable on every compact set of $\R.$ Then, \eqref{exanm4} follows since $u\in C^{0,1}_\b(\R;C_\b(\R)).$ \\
Some results avoiding this difficulty will be given in Section~\ref{secaH1}.
\end{exa}

\section{Some results of Functional Analysis}
\label{secfunana}

We recall that we shall always identify $L^2(\Omega)$ with its topological dual. Below, we recall some important results of Functional Analysis and give some new ones. Let $E$ and $F$ be locally convex Hausdorff topological vector spaces. If $E\overset{e}{\inj}F$ with dense embedding then $F^\star\overset{e^\star}{\inj}E^\star,$ where $e^\star$ is the transpose of $e:$
\begin{gather}
\label{dualtranspose}
\forall L\in F^\star, \; \forall x\in E, \; \langle e^\star(L),x\rangle_{E^\star,E}=\langle L,e(x)\rangle_{F^\star,F}.
\end{gather}
If, furthermore, $E$ is reflexive then the embedding $F^\star\overset{e^\star}{\inj}E^\star$ is dense. In most of the cases, $e$ is the identity function, so that $e^\star$ is nothing else but the restriction to $E$ of continuous linear forms on $F.$ For more details, see Trèves~\cite[Corollary~5, p.188; Corollary, p.199; Theorem~18.1, p.184]{MR2296978} and~\cite{MR4521439}. Let $A_1$ and $A_2$ be two Banach spaces such that $A_1,A_2\subset\vH$ for some Hausdorff topological vector space $\vH.$ Then $A_1\cap A_2$ and $A_1+A_2$ are Banach spaces where,
\begin{gather*}
\|a\|_{A_1\cap A_2}=\max\big\{\|a\|_{A_1},\|a\|_{A_2}\big\}	\; \text{ and } \;
\|a\|_{A_1+A_2}=\inf_{\left\{\substack{a=a_1+a_2 \hfill \\ (a_1,a_2)\in A_1\times A_2}\right.}\Big(\|a_1\|_{A_1}+\|a_2\|_{A_2}\Big).
\end{gather*}
If, in addition, $A_1\cap A_2$ is dense in both $A_1$ and $A_2$ then,
\begin{gather}
\label{dual}
\big(A_1\cap A_2\big)^\star=A_1^\star+A_2^\star \; \text{ and } \; \big(A_1+A_2\big)^\star=A_1^\star\cap A_2^\star.
\end{gather}
See, for instance, Bergh and Löfström~\cite{MR0482275} (Lemma~2.3.1 and Theorem~2.7.1). Let $Y$ be a Banach space such that $\Dr(\Omega)\inj Y$ with dense embedding. Then,
\begin{gather}
\label{injDp}
L^1_\loc\big([0,\infty);Y^\star\big)\inj\Dr^\p\big((0,\infty)\times\Omega\big).
\end{gather}
See, for instance, Droniou~\cite[Lemme~2.6.1, p.58]{droniou}. Let $I$ be an interval, let $X$ be a Banach space and let $p\in[1,\infty).$ If $X$ is separable then so is $L^p(I,X),$ and if $X$ is reflexive then,
\begin{gather*}
L^p(I;X)^\star\cong L^{p^\p}(I;X^\star).
\end{gather*}
See, for instance, Droniou~\cite{droniou} (Corollaire~1.3.2, p.13), and Edwards~\cite{MR0221256} (Theorem~8.18.3, p.590; Theorem~8.20.5, p.607). Finally, another result which will be useful is the following (Strauss~\cite[Theorem~2.1]{MR0205121}). Let $X\inj Y$ be two Banach spaces. Let $I$ be an interval and $u\in C_\w(\ovl I;Y).$ Assume that there exist $C>0$ and $N_0\subset I$ with $|N_0|=0$ such that for any $t\in I\setminus N_0,$ $u(t)\in X$ and $\|u(t)\|_X\le C.$ If $X$ is reflexive then,
\begin{gather}
\label{weakcon}
\forall t\in\ovl I, \; u(t)\in X \; \text{ and } \; u\in C_\w(\ovl I;X).
\end{gather}
Let $I$ be an interval. It is well-known that if $X$ is a Banach space which is neither reflexive, nor separable then it is not true that $L^1(I;X)$ is separable and that $L^\infty(I;X^\star)$ is the dual space of $L^1(I;X).$ In particular, we do not have a duality product $L^1(I;X)$-$L^\infty(I;X^\star)$ represented by  the natural integral. In addition, it is also well-known that the space of smooth functions is not dense in $L^\infty(I;X^\star).$ Nevertheless, we may obtain some kind of density result in ``some weak$\star$ topology", and do as if $L^\infty(I;X^\star)$ was the dual space of $L^1(I;X)$ (see, in particular, \eqref{thmdenweak4} below).

\begin{thm}[\textbf{Weak$\star$ density}]
\label{thmdenweak}
Let $\Omega\subseteq\R^N$ be an open set, let $I$ be any interval and let $X\inj X^\star$ be a Banach space. Then, for any $u\in L^1(I;X)\cap W^{1,\infty}(I;X^\star),$ there exists $(u_n)_{n\in\N}\subset\Dr(\ovl I;X)$ such that,
\begin{gather}
\label{thmdenweak1}
\forall n\in\N, \; \|u_n\|_{L^1(I;X)}\le\|u\|_{L^1(I;X)},						\\
\label{thmdenweak2}
\forall n\in\N, \; \|u_n\|_{W^{1,\infty}(I;X^\star)}\le2\|u\|_{W^{1,\infty}(I;X^\star)},	\\
\label{thmdenweak3}
u_n\xrightarrow[n\to\infty]{L^1(I;X)}u,									\\
\label{thmdenweak4}
\forall\vphi\in L^1(I;X), \; \lim_{n\to\infty}\int_I\left|\big\langle u_n^\p(t)-u^\p(t),\vphi(t)\big\rangle_{X^\star,X}\right|\d t=0.
\end{gather}
In addition, $u_n(t)\xrightarrow[n\to\infty]{X}u(t),$ for almost every $t\in I.$
\end{thm}

\begin{proof*}
With help of a continuous linear extension operator, we are brought back to the case where $I=\R.$ For the construction of such an operator, see for instance Droniou~\cite[Corollaire~2.3.1, p.48]{droniou}, Brezis~\cite[Theorem~8.6, p.209]{MR2759829}, or the Appendix in Brezis and Cazenave~\cite{bc}. Let then $u\in L^1(\R;X)\cap W^{1,\infty}(\R;X^\star).$ Let $(\rho_n)_{n\in\N}$ and $(\xi_n)_{n\in\N}$ be sequences of mollifiers and cut-off functions, respectively. The following results are standard and may be found, for instance, in Droniou~\cite{droniou} (for vector-valued functions), and Brezis~\cite{MR2759829} (for real-valued functions, but the proofs can be easily adapted for vector-valued functions). Set for any $\ell,n\in\N,$ $u_{\ell,n}=\rho_\ell\star(\xi_n u)\in\Dr(\R;X).$ By Young's inequality, we have that for any $\ell,n\in\N,$ $\|u_{\ell,n}\|_{L^1(\R;X)}\le\|u\|_{L^1(\R;X)},$ $\|u_{\ell,n}\|_{L^\infty(\R;X^\star)}\le\|u\|_{L^\infty(\R;X^\star)},$
\begin{gather}
\label{demthmdenweak1}
\eps_{\ell,n}\eqdef\|\rho_\ell\star(\xi_n u)-\xi_nu\|_{L^1(\R;X)}\xrightarrow{\ell\to\infty}0,				\\
\label{demthmdenweak2}
\eps_{\ell,n}^\p\eqdef\|\rho_\ell\star(\xi_n u^\p)-\xi_nu^\p\|_{L^1(\R;X^\star)}\xrightarrow{\ell\to\infty}0,		\\
\label{demthmdenweak3}
u_{\ell,n}^\p=\rho_\ell\star(\xi_n^\p u)+\rho_\ell\star(\xi_n u^\p),									\\
\label{demthmdenweak4}
\|\langle\rho_\ell\star(\xi_n^\p u),\vphi\rangle_{X^\star,X}\|_{L^1(\R)}\le\frac{C}n\|u\|_{L^\infty(\R;X^\star)}\|\vphi\|_{L^1(\R;X)},
\end{gather}
for any $\vphi\in L^1(\R;X).$ By \eqref{demthmdenweak1}--\eqref{demthmdenweak3}, and renumbering the sequences if necessary, we may find an increasing sequence $(\ell_n)_{n\in\N}\subset\N$ such that $\|u^\p_{\ell_n,n}\|_{L^\infty(\R;X^\star)}\le2\|u\|_{W^{1,\infty}(\R;X^\star)},$
\begin{gather}
\label{demthmdenweak5}
\eps_{\ell_n,n}+\eps_{\ell_n,n}^\p\xrightarrow{n\to\infty}0,	\\
\nonumber
\text{for almost every } t\in\R, \;\; \rho_{\ell_n}\star(\xi_n u^\p)(t)-\xi_nu^\p(t)\xrightarrow[n\to\infty]{X^\star}0.
\end{gather}
By Young's inequality, we have for any $\vphi\in L^1(\R;X)$ and $\ell,n\in\N,$
\begin{gather*}
\left|\langle\rho_\ell\star(\xi_n u^\p)-\xi_nu^\p,\vphi\rangle_{X^\star,X}\right|\le2\|u^\p\|_{L^\infty(\R,X^\star)}\|\vphi\|_X\in L^1(\R;\R),
\end{gather*}
almost everywhere in $\R.$ It follows from the dominated convergence Theorem that,
\begin{gather}
\label{demthmdenweak6}
\forall\vphi\in L^1(\R;X), \; \left\|\langle\rho_{\ell_n}\star(\xi_n u^\p)-\xi_n u^\p,\vphi\rangle_{X^\star,X}\right\|_{L^1(\R)}\xrightarrow{n\to\infty}0.
\end{gather}
Finally, still by the Lebesgue Theorem, we easily obtain that,
\begin{gather}
\label{demthmdenweak7}
\|\xi_n u-u\|_{L^1(\R;X)}\xrightarrow{n\to\infty}0,		\\
\label{demthmdenweak8}
\forall\vphi\in L^1(\R;X), \; \|\langle\xi_nu^\p-u^\p,\vphi\rangle_{X^\star,X}\|_{L^1(\R)}\xrightarrow{n\to\infty}0.
\end{gather}
Let for any $n\in\N,$ $u_n=u_{\ell_n,n}.$ Putting together \eqref{demthmdenweak5} and \eqref{demthmdenweak7}, we get \eqref{thmdenweak3}. Using \eqref{demthmdenweak3} and putting together \eqref{demthmdenweak4}, \eqref{demthmdenweak5} and \eqref{demthmdenweak8}, we get \eqref{thmdenweak4}. Finally, by \eqref{thmdenweak3}, there exists a subsequence, that we still denote by $(u_n)_{n\in\N},$ such for almost every $t\in\R,$ $u_n(t)\xrightarrow[n\to\infty]{X}u(t),$ This concludes the proof of the theorem.
\medskip
\end{proof*}

\noindent
The proof of the following result is any easy adaptation of that of Theorem~\ref{thmdenweak}, and the details are left to the reader.

\begin{thm}[\textbf{Weak$\star$ approximation}]
\label{thmdenweakbis}
Let $\Omega\subseteq\R^N$ be an open set, let $I$ be any interval and let $X\inj X^\star$ be a Banach space. Then, for any $u\in L^\infty(I;X)\cap W^{1,1}(I;X^\star),$ there exists $(u_n)_{n\in\N}\subset\Dr(\ovl I;X)$ such that,
\begin{gather}
\label{thmdenweakbis1}
\forall n\in\N, \; \|u_n\|_{W^{1,1}(I;X^\star)}\le3\|u\|_{W^{1,1}(I;X^\star)},	\\
\label{thmdenweakbis2}
\forall n\in\N, \; \|u_n\|_{L^\infty(I;X)}\le\|u\|_{L^\infty(I;X)},				\\
\label{thmdenweakbis3}
u_n\xrightarrow[n\to\infty]{W^{1,1}(I;X^\star)}u,						\\
\label{thmdenweakbis4}
\forall\vphi\in L^1(I;X^\star), \; \lim_{n\to\infty}\int_I\left|\big\langle\vphi(t),u_n(t)-u(t)\big\rangle_{X^\star,X}\right|\d t=0.
\end{gather}
In addition, $u_n(t)\xrightarrow[n\to\infty]{X^\star}u(t)$ and $u_n^\p(t)\xrightarrow[n\to\infty]{X^\star}u^\p(t),$ for almost every $t\in I.$
\end{thm}

\noindent
If $X\inj L^2(\Omega)$ with dense embedding, and if $1<p<\infty$ then for any interval $I,$ $L^p(I;X)\cap W^{1,p^\p}(I;X^\star)\inj C_\b(\ovl I;L^2(\Omega))$ (\cite[Lemma~A.4]{MR4098330}). We need and extend this to the case $p\in\{1,\infty\}.$

\begin{thm}
\label{thmembL2}
Let $\Omega\subseteq\R^N$ be an open set, let $I$ be any interval and let $X$ be a Banach space such that $X\inj L^2(\Omega)$ with dense embedding. We have the following results.
\begin{enumerate}[$1)$]
\item
\label{L1L2}
$L^1(I;X)\cap W^{1,\infty}(I;X^\star)\inj C_\b(\ovl I;L^2(\Omega)).$
\item
\label{LinfL2}
$L^\infty(I;X)\cap W^{1,1}(I;X^\star)\inj C_\b(\ovl I;L^2(\Omega)).$
\item
\label{mass}
If $u\in L^1(I;X)\cap W^{1,\infty}(I;X^\star)$ or if $u\in L^\infty(I;X)\cap W^{1,1}(I;X^\star)$ then the mapping
\begin{gather}
\label{W11}
t\longmapsto\frac12\|u(t)\|_{L^2(\Omega)}^2 \text{ belongs to } W^{1,1}(I;\R),
\end{gather}
and we have,
\begin{gather}
\label{mass'}
\frac12\frac\d{\d t}\|u(t)\|_{L^2(\Omega)}^2=\big\langle u^\p(t),u(t)\big\rangle_{X^\star,X},
\end{gather}
for almost every $t\in I.$
\end{enumerate}
\end{thm}

\begin{proof*}
By the dense embedding $X\inj L^2(\Omega),$ we have by \eqref{dualtranspose} that,
\begin{gather}
\label{demthmembL21}
X\inj L^2(\Omega)\inj X^\star,									\\
\label{demthmembL22}
\forall v\in X, \; \|v\|_{L^2(\Omega)}^2=\langle v,v\rangle_{X^\star,X}.
\end{gather}
We split the proof into four steps. \\
\textbf{Step~1:} Proof of the statement \ref{LinfL2}).
\\
Let $u\in L^\infty(I;X)\cap W^{1,1}(I;X^\star).$ By the embedding $W^{1,1}(I;X^\star)\inj C_\b(\ovl I;X^\star),$ we get by \eqref{demthmembL21}, \eqref{demthmembL22}, and \eqref{weakcon} that $u\in C_\b(\ovl I;L^2(\Omega)),$ and that there exists $C>0$ such that
\begin{gather*}
\|u(t)\|_{L^2(\Omega)}^2\le C\|u\|_{L^\infty(I;X)}\|u\|_{W^{1,1}(I,X^\star)}\le C\left(\|u\|_{L^\infty(I;X)}+\|u\|_{W^{1,1}(I,X^\star)}\right)^2,
\end{gather*}
for any $t\in\ovl I.$ Hence, $L^\infty(I;X)\cap W^{1,1}(I;X^\star)\inj C_\b(\ovl I;L^2(\Omega)).$
\\
\textbf{Step~2:} Proof of the statement \ref{mass}).
\\
Let $u\in L^1(I;X)\cap W^{1,\infty}(I;X^\star)$ $(u\in L^\infty(I;X)\cap W^{1,1}(I;X^\star),$ respectively). We claim that $\|u(\:.\:)\|_{L^2(\Omega)}\in C_\b(\ovl I;\R),$
\begin{gather}
\label{demthmembL23}
\frac12\|u(t)\|_{L^2(\Omega)}^2=\frac12\|u(s)\|_{L^2(\Omega)}^2+\int_s^t\big\langle u^\p(\sigma),u(\sigma)\big\rangle_{X^\star,X}\d\sigma,
\end{gather}
for any $t,s\in\ovl I,$ and that \eqref{W11}--\eqref{mass'} hold true. \\
For such a $u,$ let $(u_n)_{n\in\N}\subset\Dr(\ovl I;X)$ be given by Theorem~\ref{thmdenweak} (Theorem~\ref{thmdenweakbis}, respectively). Let $n\in\N.$ We have that $\|u_n(\:.\:)\|_{L^2(\Omega)}^2\in C^1(\ovl I;\R),$ and for any $t,s\in\ovl I,$
\begin{gather}
\label{demthmembL24}
\frac12\|u_n(t)\|_{L^2(\Omega)}^2=\frac12\|u_n(s)\|_{L^2(\Omega)}^2+\int_s^t\big\langle u_n^\p(\sigma),u_n(\sigma)\big\rangle_{X^\star,X}\d\sigma.
\end{gather}
By \eqref{thmdenweak2}--\eqref{thmdenweak4} (\eqref{thmdenweakbis2}--\eqref{thmdenweakbis4}, respectively), we have that,
\begin{gather}
\label{demthmembL25}
\lim_{n\to\infty}\int_I\big\langle u_n^\p(\sigma),u_n(\sigma)\big\rangle_{X^\star,X}\d\sigma
=\int_I\big\langle u^\p(\sigma),u(\sigma)\big\rangle_{X^\star,X}\d\sigma.
\end{gather}
If $u\in L^1(I;X)\cap W^{1,\infty}(I;X^\star)$ then, since $X\inj L^2(\Omega),$ we have by Theorem~\ref{thmdenweak}, that
\begin{gather}
\label{demthmembL26}
\text{for almost everywhere } t\in I, \; u_n(t)\xrightarrow[n\to\infty]{L^2(\Omega)}u(t),
\end{gather}
while if $u\in L^\infty(I;X)\cap W^{1,1}(I;X^\star),$ then \eqref{demthmembL26} comes from \eqref{demthmembL22} and Theorem~\ref{thmdenweakbis}. Passing to the limit in \eqref{demthmembL24}, it follows from \eqref{demthmembL25} and \eqref{demthmembL26} that \eqref{demthmembL23} holds true for almost every $t,s\in I.$ Now, since the mapping, $\sigma\longmapsto\big\langle u^\p(\sigma),u(\sigma)\big\rangle_{X^\star,X}$ belongs to $L^1(I;\R),$ it follows that \eqref{W11}--\eqref{mass'} come from \eqref{demthmembL23}. Finally, by \eqref{W11}--\eqref{mass'} and the embedding $W^{1,1}(I;\R)\inj C_\b(\ovl I;\R),$ we deduce that $\|u(\:.\:)\|_{L^2(\Omega)}\in C_\b(\ovl I;\R),$ and that \eqref{demthmembL23} holds true for any $t,s\in\ovl I.$
\\
\textbf{Step~3:} If $u\in L^1(I;X)\cap W^{1,\infty}(I;X^\star)$ then $u\in C_\b(\ovl I;L^2(\Omega)).$
\\
Let $u\in L^1(I;X)\cap W^{1,\infty}(I;X^\star).$ By the embedding $W^{1,\infty}(I;X^\star)\inj C_\b(\ovl I;X^\star),$ we have that $u\in C_\b(\ovl I;X^\star),$ and by \eqref{demthmembL23}, we have that $\vsup_{t\in I}\|u(t)\|_{L^2(\Omega)}<\infty.$ It then follows from \eqref{weakcon} that for any $t\in\ovl I,$ $u(t)\in L^2(\Omega),$ and $u\in C_\w(\ovl I;L^2(\Omega)).$ Let $t\in\ovl I$ and $(t_n)_{n\in\N}\subset\ovl I$ converging toward $t.$ Then by weak continuity,
\begin{gather*}
u(t_n)\underset{n\to\infty}{\overset{L^2(\Omega)_\w}{-\!\!\!-\!\!\!-\!\!\!-\!\!\!-\!\!\!\weak}}u(t).
\end{gather*}
By Step~2, we have that,
\begin{gather*}
\lim_{n\to\infty}\|u(t_n)\|_{L^2(\Omega)}=\|u(t)\|_{L^2(\Omega)}.
\end{gather*}
The space $L^2(\Omega)$ being uniformly convex, we deduce that $u(t_n)\xrightarrow[n\to\infty]{L^2(\Omega)}u(t).$ Hence, $u\in C_\b(\ovl I;L^2(\Omega)).$
\\
\textbf{Step~4:} Proof of the statement \ref{L1L2}).
\\
Let $u\in L^1(I;X)\cap W^{1,\infty}(I;X^\star).$ By Step~3, $u\in C_\b(\ovl I;L^2(\Omega)).$ Let us show the continuous embedding. It follows from \eqref{demthmembL22}, \eqref{demthmembL23}, and Hölder's and Young's inequalities that,
\begin{gather*}
\|u(t)\|_{L^2(\Omega)}^2\le\|u(s)\|_X\|u(s)\|_{X^\star}+\|u\|_{L^1(I;X)}^2+\|u^\p\|_{L^\infty(I;X^\star)}^2,
\end{gather*}
for any $t,s \in\ovl I.$ Let $(I_n)_{n\in\N}\subset I$ be a sequence of bounded intervals such that $\bigcup_{n\in\N} I_n=I.$ Integrating in $s$ and applying, once more time, Hölder's and Young's inequalities, we get,
\begin{gather*}
|I_n|\,\|u\|_{C_\b(\ovl{I_n};L^2)}^2\le(1+|I_n|)\left(\|u\|_{L^1(I;X)}+\|u\|_{W^{1,\infty}(I;X^\star)}\right)^2,
\end{gather*}
for any $n\in\N.$ Dividing by $|I_n|,$ taking the square root and letting $n\nearrow\infty,$ we arrive at,
\begin{gather*}
\|u\|_{C_\b(\ovl I;L^2)}\le(1+|I|^{-\frac12})\left(\|u\|_{L^1(I;X)}+\|u\|_{W^{1,\infty}(I;X^\star)}\right),
\end{gather*}
with the convention that $|I|^{-\frac12}=0,$ if $|I|=\infty.$ The theorem is proved.
\medskip
\end{proof*}

\noindent
Before to state some consequences we need to introduce some definitions, and auxiliary results.

\begin{defi}
\label{defRNP}
Let $\Omega\subseteq\R^N$ be an open subset. Let $Y$ be a Banach space. We shall say that a family $(Y_n)_{n\in\N_0}$ of Banach spaces is a \textit{$Y$-approximating sequence of RNP-spaces} if it satisfies the following properties.
\begin{enumerate}
\item
\label{Y1}
For any $n\in\N,$ $Y_n$ is separable and reflexive.
\item
\label{Y2}
For any $n\in\N_0,$ $Y_n\inj Y_{n+1}\inj Y,$ and for any $f\in Y_n,$ $\|f\|_Y\le\|f\|_{Y_n}.$ Moreover, each embedding is dense.
\item
\label{Y3}
For any $n\in\N_0,$ $\Dr(\Omega)\inj Y_n,$ with dense embedding.
\item
\label{Y4}
For any $f\in Y_0,$ $\vlim_{n\to\infty}\|f\|_{Y_n}=\|f\|_Y.$ 
\end{enumerate}
\end{defi}

\begin{rmk}
\label{rmkdefRNP}
Assume that $Y$ is a Banach space which admits an approximating sequence $(Y_n)_{n\in\N_0}$ of RNP-spaces. If there exists $\ell_0\in\N$ such that for any $n\ge\ell_0,$ $Y_n=Y_{\ell_0},$ then it follows from Definition~\ref{defRNP} that $Y$ and $Y_{\ell_0}$ are two Banach spaces with the same norm, and that $Y_{\ell_0}$ is dense in $Y.$ Therefore, $Y_{\ell_0}=Y.$ As a consequence, if $Y$ is not separable or not reflexive then, renumbering $(Y_n)_{n\in\N_0}$ if necessary, we have that for any $n\in\N_0,$ $Y_n\subsetneq Y_{n+1}\subsetneq Y.$
\end{rmk}

\begin{lem}
\label{lemY*}
Let $\Omega\subseteq\R^N$ be an open subset. Let $Y$ be a Banach space and let $(Y_n)_{n\in\N_0}$ be a $Y$-approximating sequence of RNP-spaces. Then $(Y_n^\star)_{n\in\N_0}$ satisfies the following properties.
\begin{enumerate}
\item
\label{Y*1}
For any $n\in\N,$ $Y_n^\star$ is separable and reflexive.
\item
\label{Y*2}
For any $n\in\N_0,$ $Y^\star\inj Y_{n+1}^\star\inj Y_n^\star,$ and for any $f\in Y^\star,$ $\|f\|_{Y_n^\star}\le\|f\|_{Y^\star}.$ Moreover, if $n\in\N$ then each embedding is dense.
\item
\label{Y*3}
For any $n\in\N_0,$ $Y_n^\star\inj\Dr^\p(\Omega),$ with dense embedding.
\item
\label{Y*4}
For any $f\in Y^\star,$ $\vlim_{n\to\infty}\|f\|_{Y_n^\star}=\|f\|_{Y^\star}.$ 
\end{enumerate}
\end{lem}

\begin{proof*}
It is well-known that the dual space of a reflexive and separable Banach space is also a reflexive and separable Banach space. Therefore, Property~\ref{Y*1} holds true. The rest of the lemma follows easily by duality, Definition~\ref{defRNP}, and the results about Functional Analysis we recalled at the beginning of this section.
\medskip
\end{proof*}

\begin{cor}
\label{corZY}
Let $\Omega\subseteq\R^N$ be an open subset. Let $Y$ be Banach space and let $(Y_n)_{n\in\N_0}$ be a $Y$-approximating sequence of RNP-spaces. Finally, let $Z$ be a separable and reflexive Banach space such that $\Dr(\Omega)\inj Z,$ with dense embedding. Then, $(Z\cap Y_n)_{n\in\N_0}$ is a $Z\cap Y$-approximating sequence of RNP-spaces.
\end{cor}

\begin{proof*}
Let $n\in\N.$ Since $\Dr(\Omega)\inj Z$ and $\Dr(\Omega)\inj Y_n$ with dense embeddings, it follows that convergences in the weak topologies $Z_\w$ and $Y_{n,\w}$ imply convergence in $\Dr^\p(\Omega).$ Therefore, we easily obtain from \eqref{dual} and the Eberlein-\v{S}mulian Theorem that $Z\cap Y_n$ is reflexive. For $n\in\N,$ let $T:Z\cap Y_n\tends Z\times Y_n$ be defined by $T(u)=(u,u).$ We recall that for any $(u,v)\in Z\times Y_n,$ $\|(u,v)\|_{Z\times Y_n}=\max\big\{\|u\|_Z,\|v\|_{Y_n}\big\}.$ It is clear that $T$ is an isometry and that $Z\times Y_n$ is separable. It follows that $T(Z\cap Y_n)$ is separable (Brezis~\cite[Proposition~3.25, p.73]{MR2759829}), and so is $Z\cap Y_n.$ Then, Property~\ref{Y1} of Definition~\ref{defRNP} is satisfied. Let $n\in\N_0.$ It is clear that $\Dr(\Omega)\inj Z\cap Y_n.$ Let us show that this embedding is dense. Let $T\in Z^\star+Y_n^\star$ be such that for any $\vphi\in\Dr(\Omega),$ $\langle T,\vphi\rangle_{Z^\star+Y_n^\star,Z\cap Y_n}=0.$ We write $T=T_1+T_2,$ where $(T_1,T_2)\in Z^\star\times Y_n^\star.$ Using the dense embeddings $Z\cap Y_n\inj Z,$ $Z\cap Y_n\inj Y_n,$ $\Dr(\Omega)\inj Z$ and $\Dr(\Omega)\inj Y_n,$ it follows from \eqref{dualtranspose} that for any $\vphi\in\Dr(\Omega),$
\begin{align*}
 0=	&	\; \langle T,\vphi\rangle_{Z^\star+Y_n^\star,Z\cap Y_n}=\langle T_1,\vphi\rangle_{Z^\star,Z}+\langle T_2,\vphi\rangle_{Y_n^\star,Y_n}	\\
   =	&	\; \langle T_1,\vphi\rangle_{\Dr^\p(\Omega),\Dr(\Omega)}+\langle T_2,\vphi\rangle_{\Dr^\p(\Omega),\Dr(\Omega)}
			=\langle T_1+T_2,\vphi\rangle_{\Dr^\p(\Omega),\Dr(\Omega)}														\\
   =	&	\; \langle T,\vphi\rangle_{\Dr^\p(\Omega),\Dr(\Omega)}.
\end{align*}
Therefore, $T=0$ in $\Dr^\p(\Omega),$ hence in $Z^\star+Y_n^\star.$ The rest of the proof is obvious.
\medskip
\end{proof*}

\begin{lem}
\label{lemY}
Let $\Omega\subseteq\R^N$ be any open subset. Then $L^1(\Omega)$ admits an approximating sequence $(Y_n)_{n\in\N_0}$ of RNP-spaces. In addition, $Y_0$ may be chosen separable.
\end{lem}

\begin{proof*}
Let $\Omega\subseteq\R^N$ be an open subset. For convenience, we introduce some notations. Let $n\in\N,$ $\eps_n=\frac1n,$ $\Omega_n=\{x\in\Omega;|x|\le n\},$ and let $\omega_{N-1}$ be the area of the unit sphere $\vsS_{N-1}$ of $\R^N$ (with the convention that $\omega_0=2).$ We define,
\begin{align*}
&	Y_n=\left\{f\in L^{1+\eps_n}(\Omega);|f|^{1+\eps_n}|\:.\:|^{\eps_n(N+\eps_n)}\1_{\Omega_n^\co}\in L^1(\Omega)\right\},	\\
&	Y_0=\left\{f\in Y_1;f\,|\:.\:|^{\eps_1(N+\eps_1)}\in L^1(\Omega)\right\},
\end{align*}
whose norms are,
\begin{align*}
&	\|f\|_{Y_n}=(2\omega_{N-1}n^N)^\frac1{n+1 }\left(\int_{\Omega_n}|f|^{1+\eps_n}\d x
			+\int_{\Omega_n^\co}|f|^{1+\eps_n}|x|^{\eps_n(N+\eps_n)}\d x\right)^\frac1{1+\eps_n},	\\
&	\|f\|_{Y_0}=\max\left\{\|f\|_{Y_1},\int_\Omega|f||x|^{\eps_1(N+\eps_1)}\d x\right\}.
\end{align*}
Note that $Y_n\inj L^{1+\eps_n}(\Omega).$ Let $\lambda$ be the Lebesgue measure on the Lebesgue sets $\Br(\Omega)$ of $\Omega,$ let
\begin{align*}
&	g_n(x)=(2\omega_{N-1}n^N)^\frac1n\left(\1_{\Omega_n}(x)+|x|^{\eps_n(N+\eps_n)}\1_{\Omega_n^\co}(x)\right),	\\
&	g_0(x)=|x|^{\eps_1(N+\eps_1)},
\end{align*}
for any $x\in\Omega,$ and let $\nu_n$ be the density measure defined on $\Br(\Omega)$ by
\begin{gather*}
\frac{\d\nu_n}{\d\lambda}=g_n \; \text{ and } \; \frac{\d\nu_0}{\d\lambda}=g_0.
\end{gather*}
It follows that $Y_n=L^{1+\eps_n}(\Omega,\Br(\Omega),\nu_n),$ and $Y_0=Y_1\cap L^1(\Omega,\Br(\Omega),\nu_0).$ Classical results on measure and integration theory give that $Y_n$ is a separable and reflexive Banach space, and that $Y_0$ is a separable Banach space (for the proof of the separability of $Y_0,$ proceed as in Corollary~\ref{corZY}). Therefore, $(Y_n)_{n\in\N_0}$ satisfies Property~\ref{Y1} of Definition~\ref{defRNP}. Let $n\in\N.$ It is also obvious that $\Dr(\Omega)\inj Y_0\inj Y_n.$ Let us establish the density. Let $f\in Y_n.$ Set $a=\eps_n(N+\eps_n)$ and $b=1+\eps_n.$ Let $(\rho_j)_{j\in\N}$ be a sequence of mollifiers. Let $O_\ell=\left\{x\in\Omega_\ell;\dist(x,\Omega^\co)>\frac2\ell\right\},$ $\ell\in\N.$ Denote by $\wt f$ the extension of $f$ by $0$ outside of $\Omega.$ Finally, for $j,\ell\in\N,$ let $\wt{\vphi_j^\ell}=\rho_j\star(\wt f\1_{O_\ell}).$ For any $\ell\in\N,$ and $j>\ell,$ we have that
\begin{align}
\label{prooflemY1}
2^{1-b}\int_{\Omega_n^\co}|x|^a|\wt{\vphi_j^\ell}-\wt f|^b\d x	&	\;
\le(2\ell)^a\int_{\Omega_n^\co}|\wt{\vphi_j^\ell}-\wt f\1_{O_\ell}|^b\d x+\int_{\Omega_n^\co}|x|^a|\wt f\1_{O_\ell}-\wt f|^b\d x,	\\
\label{prooflemY2}
2^{1-b}\int_{\Omega_n}|\wt{\vphi_j^\ell}-\wt f|^b\d x	&	\;
\le\int_{\Omega_n}|\wt{\vphi_j^\ell}-\wt f\1_{O_\ell}|^b\d x+\int_{\Omega_n}|\wt f\1_{O_\ell}-\wt f|^b\d x.
\end{align}
By the Lebesgue Theorem, the last integrals in \eqref{prooflemY1} and \eqref{prooflemY2} go to $0,$ as $\ell\to\infty,$ while the first of the right members go to $0,$ as $j\to\infty,$ for each $\ell\in\N$ (by the classical results about truncation and regularization). It follows that there exists an increasing sequence $(j_\ell)_{\ell\in\N}$ which goes to $\infty,$ as $\ell\tends\infty,$ such that
\begin{gather*}
\lim_{\ell\to\infty}\int_{\Omega_n}|\wt{\vphi_{j_\ell}^\ell}-\wt f|^b\d x=\lim_{\ell\to\infty}\int_{\Omega_n^\co}|x|^a|\wt{\vphi_{j_\ell}^\ell}-\wt f|^b\d x=0.
\end{gather*}
Let for any $\ell\in\N,$ $\vphi_\ell=\wt{\vphi_{j_\ell}^\ell}_{|\Omega}\in\Dr(\Omega).$ It follows that $(\vphi_\ell)_{\ell\in\N}$ answers to the problem. A trivial adaptation of the proof gives the density of $\Dr(\Omega)$ in $Y_0.$ Therefore, $(Y_n)_{n\in\N_0}$ satisfies Property~\ref{Y3} of Definition~\ref{defRNP}. Now, let us show that $(Y_n)_{n\in\N_0}$ satisfies Property~\ref{Y2}. The density result comes from the density of $\Dr(\Omega)$ in $L^1(\Omega)$ and Property~\ref{Y3}. Let $n\in\N$ (the case $n=0$ is immediate). Let $f\in Y_n.$ We have,
\begin{align*}
	&	\; (2\omega_{N-1}(n+1)^N)^{-\frac1{n+1}}\|f\|_{Y_{n+1}}^{1+\eps_{n+1}}											\\
   =	&	\; \int_{\Omega_{n+1}}|f|^{1+\eps_{n+1}}\d x+\int_{\Omega_{n+1}^\co}|f|^{1+\eps_{n+1}}|x|^{\eps_{n+1}(N+\eps_{n+1})}\d x	\\
  \le	&	\; \int_{\Omega_n}|f|^{1+\eps_{n+1}}\d x+\int_{\Omega_n^\co}|f|^{1+\eps_{n+1}}|x|^{\eps_{n+1}(N+\eps_{n+1})}\d x.
\end{align*}
By Hölder's inequality,
\begin{gather*}
\int_{\Omega_n}|f|^{1+\eps_{n+1}}\d x\le|B(0,n)|^\frac1{(n+1)^2}\left(\int_{\Omega_n}|f|^{1+\eps_n}\d x\right)^\frac{1+\eps_{n+1}}{1+\eps_n},
\end{gather*}
and
\begin{align*}
	&	\; \int_{\Omega_n^\co}|f|^{1+\eps_{n+1}}|x|^{\eps_{n+1}(N+\eps_{n+1})}\d x									\\
   =	&	\; \int_{\Omega_n^\co}|x|^{\eps_{n+1}(N+\eps_{n+1})-\eps_n(N+\eps_n)\frac{1+\eps_{n+1}}{1+\eps_n}}
			\left(|f|^{1+\eps_{n+1}}|x|^{\eps_n(N+\eps_n)\frac{1+\eps_{n+1}}{1+\eps_n}}\right)\d x						\\
  \le	&	\; \left(\int_{\Omega_n^\co}|x|^{-(N+\eps_n\eps_{n+1}+\eps_n+\eps_{n+1})}\d x\right)^\frac{\eps_n-\eps_{n+1}}{1+\eps_n}
			\left(\int_{\Omega_n^\co}|f|^{1+\eps_n}|x|^{\eps_n(N+\eps_n)}\d x\right)^\frac{1+\eps_{n+1}}{1+\eps_n}			\\
  \le	&	\; C(N,n)\left(\int_{\Omega_n^\co}|f|^{1+\eps_n}|x|^{\eps_n(N+\eps_n)}\d x\right)^\frac{1+\eps_{n+1}}{1+\eps_n}.
\end{align*}
Gathering together the above estimates, we get that, $\|f\|_{Y_{n+1}}\le C(N,n)\|f\|_{Y_n}.$ Now, let $f\in Y_n.$ Let us show that $f\in L^1(\Omega)$ and $\|f\|_{L^1(\Omega)}\le\|f\|_{Y_n}.$ Using twice Hölder's inequality, we have that,
\begin{align*}
	&	\; \|f\|_{L^1(\Omega)}		\\
  \le	&	\; |B(0,n)|^\frac{\eps_n}{1+\eps_n}\left(\int_{\Omega_n}|f|^{1+\eps_n}\d x\right)^\frac1{1+\eps_n}
			+\int_{\Omega_n^\co}|x|^{-\frac{\eps_n(N+\eps_n)}{1+\eps_n}}\left(|f||x|^\frac{\eps_n(N+\eps_n)}{1+\eps_n}\right)\d x		\\
  \le	&	\; |B(0,n)|^\frac{\eps_n}{1+\eps_n}\left(\int_{\Omega_n}|f|^{1+\eps_n}\d x\right)^\frac1{1+\eps_n}						\\
	&		\qquad	+\omega_{N-1}^\frac{\eps_n}{1+\eps_n}\left(\int_n^\infty r^{-(1+\eps_n)}\d r\right)^\frac{\eps_n}{1+\eps_n}
			\left(\int_{\Omega_n^\co}|f|^{1+\eps_n}|x|^{\eps_n(N+\eps_n)}\d x\right)^\frac1{1+\eps_n}							\\
  \le	&	\; (\omega_{N-1}n^N)^\frac1{n+1}\left(\int_{\Omega_n}|f|^{1+\eps_n}\d x\right)^\frac1{1+\eps_n}							\\
	&		\qquad	+\omega_{N-1}^\frac1{n+1}n^{-\frac1n}
			\left(\int_{\Omega_n^\co}|f|^{1+\eps_n}|x|^{\eps_n(N+\eps_n)}\d x\right)^\frac1{1+\eps_n}							\\
  \le	&	\;  (\omega_{N-1}n^N)^\frac1{n+1}\left(\left(\int_{\Omega_n}|f|^{1+\eps_n}\d x\right)^\frac1{1+\eps_n}
			+\left(\int_{\Omega_n^\co}|f|^{1+\eps_n}|x|^{\eps_n(N+\eps_n)}\d x\right)^\frac1{1+\eps_n}\right)					\\
  \le	&	\; \|f\|_{Y_n},
\end{align*}
since for any $a,b\ge0$ and $0<\alpha<1,$ $a^\alpha+b^\alpha\le2^{1-\alpha}(a+b)^\alpha.$ Hence Property~\ref{Y3} is satisfied. Now, to prove Property~\ref{Y4}, it is sufficient to see that for any $f\in Y_0,$ and $n\in\N,$
\begin{align*}
	&	\; |f|^{1+\eps_n}\1_{\Omega_n}\le|f|^{1+\eps_1}\1_{\{|f|>1\}}+|f|\1_{\{|f|\le1\}}\in L^1(\Omega),	\\
	&	\; |f|^{1+\eps_n}|\:.\:|^{\eps_n(N+\eps_n)}\1_{\Omega_n^\co}					\\
  \le	&	\; |f|^{1+\eps_1}|\:.\:|^{\eps_1(N+\eps_1)}\1_{\Omega_1^\co}\1_{\{|f|>1\}}
			+|f||\:.\:|^{\eps_1(N+\eps_1)}\1_{\Omega_1^\co}\1_{\{|f|\le1\}}\in L^1(\Omega),
\end{align*}
and to apply the dominated convergence Theorem. This ends the proof of the lemma.
\medskip
\end{proof*}

\section{Proofs of the theorems on uniqueness and existence of strong and weak solutions}
\label{proofexi}

\begin{vproof}{of Proposition~\ref{propsolL2}.}
Let $(u,U,f)$ be a weak solution and let $(u_n,U_n,f)_{n\in\N}$ satisfying \eqref{fn}--\eqref{U}. Let $(Y_\ell)_{\ell\in\N_0}$ be any $L^1$-approximating sequence of RNP-spaces, which exists by Lemma~\ref{lemY}. It follows from Corollary~\ref{corZY} that $(H^2_0\cap Y_\ell)_{\ell\in\N_0}$ is an $H^2_0\cap L^1$-approximating sequence of RNP-spaces. By \eqref{cv}, \cite[Lemma~4.2]{MR4340780}, and the diagonal procedure, we have that (up to a subsequence),
\begin{align}
\label{dempropsolL21}
&	\Delta u_n\xrightarrow[n\to\infty]{C([0,T];H^{-2}(\Omega))}\Delta u,	\\
\label{dempropsolL22}
&	Vu_n\xrightarrow[n\to\infty]{C([0,T];H^{-1}(\Omega))}Vu,			\\
\label{dempropsolL23}
&	u_n\xrightarrow[n\to\infty]{\text{a.e.\,in }(0,\infty)\times\Omega}u.
\end{align}
for any $T>0.$ Let $T>0,$ $k\in\N,$ $\Omega_k=\Omega\cap B(0,k),$ and $\omega_{T,k}=\big\{(t,x)\in(0,T)\times\Omega_k; u(t,x)\neq0\big\}.$ Finally, let $\omega=\big\{(t,x)\in(0,\infty)\times\Omega; u(t,x)\neq0\big\}.$ Note that by \eqref{U}, we have,
\begin{gather}
\label{dempropsolL24}
U_{n|\omega_{T,k}}\underset{n\to\infty}{\overset{L^2(\omega_{T,k})_\w}{-\!\!\!-\!\!\!-\!\!\!-\!\!\!-\!\!\!-\!\!\!-\!\!\!-\!\!\!\weak}}U_{|\omega_{T,k}}.
\end{gather}
On the other hand, it follows from \eqref{dempropsolL23} and the dominated convergence Theorem that,
\begin{gather}
\label{dempropsolL25}
\frac{u_n}{|u_n|}_{|\omega_{T,k}}\xrightarrow[n\to\infty]{L^2(\omega_{T,k})}\frac{u}{|u|}_{|\omega_{T,k}}.
\end{gather}
Finally, by \eqref{defsol23} and \eqref{dempropsolL23}, we have for any $n$ large enough that $U_n=\frac{u_n}{|u_n|},$ almost everywhere in $\omega_{T,k}.$ Since $T$ and $k$ are arbitrary, we then deduce from \eqref{dempropsolL24} and \eqref{dempropsolL25} that $U=\frac{u}{|u|},$ almost everywhere in $\omega.$ Therefore, $U$ satisfies \eqref{defsol22}--\eqref{defsol23}. Now, with help of \eqref{cv}, \eqref{U}, \eqref{dempropsolL21} and \eqref{dempropsolL22}, we have that $u$ satisfies \eqref{nls} in $\Dr^\p\big((0,\infty)\times\Omega\big).$ Let $\ell\in\N.$ By Definition~\ref{defsol}, $(u_n)_{n\in\N}\subset W^{1,\infty}_\loc\big([0,\infty);X^\star\big),$ where $X=H^1_0(\Omega)\cap L^1(\Omega).$ By \eqref{nls}, it follows that for any $n\in\N,$ $t\longmapsto U_n(t)$ is measurable $[0,\infty)\tends H^{-2}(\Omega)+L^\infty(\Omega).$ By Corollary~\ref{corZY}, $(U_n)_{n\in\N}$ is bounded in $L^\infty(0,\infty;H^{-2}(\Omega)+Y_\ell^\star)$ (whose each norm is bounded by $1).$ But $H^2_0(\Omega)\cap Y_\ell$ is separable and reflexive, so that $L^1(0,\infty;H^2_0(\Omega)\cap Y_\ell)$ is separable, and
\begin{gather*}
L^1(0,\infty;H^2_0(\Omega)\cap Y_\ell)^\star\cong L^\infty(0,\infty;H^{-2}(\Omega)+Y_\ell^\star)\inj\Dr^\p\big((0,\infty)\times\Omega\big).
\end{gather*}
See the beginning of Section~\ref{secfunana}. And since the limit in \eqref{U} also takes place in $\Dr^\p\big((0,\infty)\times\Omega\big),$ we deduce that
\begin{gather*}
U\in L^\infty(0,\infty;H^{-2}(\Omega)+Y_\ell^\star).
\end{gather*}
It then follows from \eqref{nls} that \eqref{propsolL21} holds true, and that \eqref{nls} makes sense in $L^1_\loc\big([0,\infty);H^{-2}(\Omega)+Y_\ell^\star\big),$ for any $\ell\in\N_0.$ In addition, \eqref{propsolL22} comes easily from \eqref{nls}, and it is clear that $u$ solves~\eqref{nls} in $ H^{-2}(\Omega)+L^\infty(\Omega),$ for almost every $t>0.$ Finally, note that if $|\Omega|<\infty$ then $H^2_0(\Omega)\cap L^1(\Omega)=H^2_0(\Omega),$ which is reflexive and separable. In this case, the above arguments work for $H^2_0(\Omega)$ in place of $H^2_0(\Omega)\cap Y_\ell.$
\medskip
\end{vproof}

\noindent
From now, we suppose Assumption~\ref{ass1}. Before proving the other results of Section~\ref{exiuni}, we recall some results of our previous papers we will need. Here and in the rest of this article, we shall use the following notations and conventions.

\noindent
For any $u\in L^2(\Omega),$ $Vu\in H^{-1}(\Omega)$ and for any $u\in H^1_0(\Omega),$ $Vu\in L^2(\Omega).$ There exists $C=C(N,\beta)>0$ such that for any $u\in H^1_0(\Omega),$
\begin{gather}
\label{lemVL2}
\|Vu\|_{L^2(\Omega)}\le C\|V\|_{L^\infty(\Omega)+L^{p_V}(\Omega)}\|u\|_{H^1_0(\Omega)}.
\end{gather}
See \cite[Lemma~4.1]{MR4340780}.

\noindent
Let $\eps\ge0.$ For any $u\in L^0(\Omega)$ and almost every $x\in\Omega,$ we define
\begin{align*}
&	g_\eps(u)(x)=(|u(x)|^2+\eps)^{-\frac12}u(x), \; \eps>0,	\\
&	g(u)(x)=g_0(u)(x)=\frac{u(x)}{|u(x)|}, \; u(x)\neq0.
\end{align*}
Now, for $\mu>0,$ let us define the operator $(A,D(A))$ on $L^2(\Omega)$ given by,
\begin{align*}
&	D(A)=\left\{u\in H^1_0(\Omega)\cap L^1(\Omega); \left(-\vi\Delta u+\mu\dfrac{u}{|u|}\right)_{|\{u\neq0\}}\in L^2(\{u\neq0\})\right\},	\medskip \\
&	C(u)=\left\{U\in\ovl B_{L^\infty}(0,1); -\vi\Delta u+\mu\,U\in L^2(\Omega) \text{ and if } u(x)\neq0, \; U(x)=\frac{u(x)}{|u(x)|}\right\},		\medskip \\
&	Au=\Big\{-\vi\Delta u-\vi Vu+\mu\,U; U\in C(u)\Big\}, \; \forall u\in D(A).
\end{align*}

\medskip
\noindent
Notice that, as in the theory of maximal monotone operators (Brezis~\cite{MR0348562}) the operator $A$ is mulitvalued, since, at least formally $A0=\ovl B_{L^\infty}(0,1)\cap L^2(\Omega).$ Note that by the inclusion $L^\infty(\Omega)\subset L^2_\loc(\Omega),$ we have $D(A)\subset H^2_\loc(\Omega).$ We split the proof of Theorem~\ref{thmstrongaH2} into several lemmas. The proofs of Corollary~\ref{corAmon} and Lemma~\ref{lemAmax} below are close to those of \cite[Corollary~5.9 and Lemma~5.12]{MR4340780}. However, they require an adaptation because we do not have that $\Delta u$ and $\frac{u}{|u|}$ belong separately in $L^2(\Omega).$ \\
Let us start by proving that the operator $A$ is monotone in $L^2(\Omega).$

\begin{lem}
\label{lemAmon}
Let $u_1,u_2\in L^1(\Omega)$ and $U_1,U_2\in\ovl{B}_{L^\infty}(0,1)$ be such that for any $j\in\{1,2\},$ $U_j=\frac{u_j}{|u_j|},$ almost everywhere where $u_j\neq0.$ Then,
\begin{gather*}
\Re\left(\int_\Omega(U_1-U_2)\big(\ovl{u_1-u_2}\big)\d x\right)\ge0.
\end{gather*}
\end{lem}

\begin{proof*}
Let for $j\in\{1,2\},$ $\omega_j=\big\{x\in\Omega; \; u_j(x)\neq0\big\}.$ We have that,
\begin{align*}
		&	\; \Re\left(\int_\Omega(U_1-U_2)\big(\ovl{u_1-u_2}\big)\d x\right)													\\
	=	&	\; \Re\left(\int_{\omega_1^\co\cap\,\omega_2}\left(U_1-\frac{u_2}{|u_2|}\right)\ovl{(-u_2)}\d x\right)
				+\Re\left(\int_{\omega_1\cap\,\omega_2^\co}\left(\frac{u_1}{|u_1|}-U_2\right)\ovl{u_1}\d x\right)							\\
		&		\qquad \; +\Re\left(\int_{\omega_1\cap\,\omega_2}\left(\frac{u_1}{|u_1|}-\frac{u_2}{|u_2|}\right)\big(\ovl{u_1-u_2}\big)\d x\right)	\\
\eqdef	&	\; I_1+I_2+I_3.
\end{align*}
Since $|U_1\ovl{u_2}|\le|u_2|$ and $|U_2\ovl{u_1}|\le|u_1|,$ we get that, $I_1\ge0$ and $I_2\ge0.$ Finally, $I_3\ge0$ by \cite[Corollary~5.5]{MR4340780}, and the lemma is proved.
\medskip
\end{proof*}

\begin{cor}
\label{corAmon}
$(A,D(A))$ is monotone on $L^2(\Omega).$
\end{cor}

\begin{proof*}
Let $(u_1,u_2)\in D(A)\times D(A)$ and $(W_1,W_2)\in Au_1\times Au_2.$ Let then $(U_1,U_2)\in C(u_1)\times C(u_2)$ be such that $-\vi\Delta u_j-\vi Vu_j+\mu\,U_j=W_j,$ for any $j\in\{1,2\}.$ Since $\Delta u_j$ or $U_j$ may not belong separately to $L^2(\Omega),$ we have to proceed in a different way than in \cite{MR4340780}. By \cite[Lemma~4.4]{MR4503241},
\begin{align*}
	&	\; (W_1-W_2,u_1-u_2)_{L^2(\Omega)}=(-\vi V(u_1-u_2),u_1-u_2)_{L^2(\Omega)}	\\
	&	\; \quad +\big(\vi\nabla(u_1-u_2),\nabla(u_1-u_2)\big)_{L^2(\Omega)}+\mu\langle U_1-U_2,u_1-u_2\rangle_{L^\infty(\Omega),L^1(\Omega)}\\
   =	&	\; \mu\,\Re\left(\int_\Omega(U_1-U_2)\big(\ovl{u_1-u_2}\big)\d x\right),
\end{align*}
and we conclude with help of Lemma~\ref{lemAmon}.
\medskip
\end{proof*}

\begin{lem}
\label{lemAmax}
$R(I+A)=L^2(\Omega).$
\end{lem}

\begin{proof*}
Let $F\in L^2(\Omega).$ By \cite[Lemma~4.3]{MR4503241}, there exist $u\in H^1_0(\Omega)\cap L^1(\Omega)$ (hence $Vu\in L^2(\Omega),$ by \eqref{lemVL2}) and a sequence $(u_{\eps_n})_{n\in\N}\subset H^1_0(\Omega)$ with $(\Delta u_{\eps_n})_{n\in\N}\subset L^2(\Omega),$ where $(\eps_n)_{n\in\N}\subset(0,\infty)$ is a decreasing sequence converging toward $0,$ satisfying the following properties: for each $n\in\N,$ $u_{\eps_n}$ is the unique solution to,
\begin{gather}
\label{lemAmax1}
-\vi\Delta u_{\eps_n}-\vi Vu_{\eps_n}+\mu g_{\eps_n}(u_{\eps_n})+u_{\eps_n}=F, \text{ in } L^2(\Omega),
\end{gather}
and
\begin{align}
\label{lemAmax2}
&	u_{\eps_n}\xrightarrow[n\to\infty]{\Dr^\p(\Omega)}u,		\\
\label{lemAmax3}
&	Vu_{\eps_n}\xrightarrow[n\to\infty]{\Dr^\p(\Omega)}Vu,	\\
\label{lemAmax4}
&	u_{\eps_n}\xrightarrow[n\to\infty]{\text{a.e.\,in }\Omega}u.
\end{align}
Now, let us denote by $\omega=\{x\in\Omega; u(x)\neq0\},$ and for $k\in\N,$ let $\omega_k=\omega\cap B(0,k).$ Let $k\in\N.$ Since $\big(g_{\eps_n}(u_{\eps_n})\big)_{n\in\N}$ belongs to $\ovl B_{L^\infty}(0,1),$ there exists $U\in \ovl B_{L^\infty}(0,1)$ such that, up to a subsequence,
\begin{gather}
\label{lemAmax5}
g_{\eps_n}(u_{\eps_n})\underset{n\to\infty}{-\!\!\!-\!\!\!-\!\!\!-\!\!\!\weak}U \text{ in } L^\infty(\Omega)_{\w\star}.
\end{gather}
In addition, we have that $g_{\eps_n}(u_{\eps_n})\xrightarrow[n\to\infty]{\text{a.e.\,in }\omega}g(u),$ by \eqref{lemAmax4}. It follows from the dominated convergence Theorem that,
\begin{gather}
\label{lemAmax6}
g_{\eps_n}(u_{\eps_n})_{|\omega_k}\xrightarrow[n\to\infty]{L^1(\omega_k)}g(u)_{|\omega_k}.
\end{gather}
Let $h\in L^\infty(\Omega)\cap L^1(\Omega)$ be defined by,
\begin{gather*}
h=
\begin{cases}
g(u)-U,	&	\text{in } \omega_k,				\\
0,		&	\text{in } \Omega\setminus\omega_k.
\end{cases}
\end{gather*}
We have by Hölder's inequality that,
\begin{align*}
	&	\; \vint_{\omega_k}|g(u)-U|^2\d x																				\\
   =	&	\; \Re\vint_{\omega_k}\big(g(u)-g_{\eps_n}(u_{\eps_n})\big)\ovl h\,\d x+\Re\vint_{\omega_k}\big(g_{\eps_n}(u_{\eps_n})-U\big)\ovl h\,\d x	\\
  \le	&	\; 2\|g_{\eps_n}(u_{\eps_n})-g(u)\|_{L^1(\omega_k)}+\langle g_{\eps_n}(u_{\eps_n})-U,h\rangle_{L^\infty(\Omega),L^1(\Omega)}.
\end{align*}
With help of \eqref{lemAmax5} and \eqref{lemAmax6}, we are allowed to pass to the limit in the above to infer that for any $k\in\N,$ $U=g(u),$ almost everywhere in $\omega_k.$ Hence,
\begin{gather}
\label{lemAmax7}
U=g(u), \text{ a.e. in } \omega.
\end{gather}
Now, passing to the limit in \eqref{lemAmax1}, we get with help of \eqref{lemAmax2}, \eqref{lemAmax3}, \eqref{lemAmax5} and \eqref{lemAmax7} that
\begin{gather*}
u\in D(A)	\; \text{ and } \;	(I+A)u\ni F.
\end{gather*}
This ends the proof.
\medskip
\end{proof*}

\begin{rmk}
\label{rmkADE}
In \cite{MR4503241}, we study the sublinear problem
\begin{empheq}[left=\empheqlbrace]{align}
	\label{nls*}
	\vi\frac{\partial u}{\partial t}+\Delta u+V(x)u+ag(u)=f(t,x),	&	\text{ in } (0,\infty)\times\Omega,				\\
	\label{nlsb*}
	u_{|\partial\Omega}=0,							&	\text{ on } (0,\infty)\times\partial\Omega,	\dfrac{}{}	\\
	\label{u0*}
	u(0)= u_0,										&	\text{ in } \Omega,
\end{empheq}
where $g(u)=|u|^{-(1-m)}u,$ with $0<m<1$ (and under some suitable assumptions about $V,$ $a\in\C$ and $f).$ To prove existence of a solution with initial data in $H^1_0(\Omega)$ (see \cite[Theorem~2.10]{MR4503241}), we proceed by regularizing the nonlinearity $g(u)$ with $g_\eps^m(u)\eqdef(|u|^2+\eps)^{-\frac{1-m}2}u,$ $\eps>0.$ This proof contains a slight flaw. Indeed, taking the $L^2$-scalar product of the regularized equation,
\begin{gather*}
\vi\frac{\partial u_\eps}{\partial t}+\Delta u_\eps+V(x)u_\eps+ag_\eps^m(u_\eps)=f_\eps(t,x),	\text{ in } L^2(\Omega),
\end{gather*}
with $\vi u_\eps,$ and applying the Cauchy-Schwarz inequality, we wrongly arrive at
\begin{gather}
\label{AdemlemthmsaH11}
\frac12\frac{\d}{\d t}\|u_\eps(\sigma)\|_{L^2(\Omega)}^2+\Im(a)\|u_\eps(\sigma)\|_{L^{m+1}(\Omega)}^{m+1}
\le\|f_\eps(\sigma)\|_{L^2(\Omega)}\|u_\eps(\sigma)\|_{L^2(\Omega)}.
\end{gather}
But actually, the correct result is,
\begin{gather}
\label{AdemlemthmsaH12}
\frac12\frac{\d}{\d t}\|u_\eps(\sigma)\|_{L^2(\Omega)}^2+\Im(a)\vint_\Omega\frac{|u_\eps(\sigma,x)|^2}{(|u_\eps(\sigma,x)|^2+\eps)^\frac{1-m}2}\d x
\le\|f_\eps(\sigma)\|_{L^2(\Omega)}\|u_\eps(\sigma)\|_{L^2(\Omega)},
\end{gather}
and then \cite[Lemmas~4.6 and 4.7]{MR4503241}, and their proofs have to be modified. The statement of \cite[Lemma~4.8]{MR4503241} is totally unchanged but its proof has to be very slightly changed. \cite[Lemmas~4.6--4.8]{MR4503241} are only needed to prove \cite[Theorems~2.9 and 2.10]{MR4503241}. Once these lemmas are correctly stated and proved, the proof of these theorems is totally unchanged. The correct version follows with very closed arguments to the ones used in \cite{MR4503241}. The only important modification consists in to replace formula (4.32) of \cite[Lemma~4.6]{MR4503241} by properties
\begin{gather}
\label{AlemthmsaH15}
\begin{cases}
(u_\eps)_{\eps>0} \text{ is bounded in }
	L^\infty_\loc\big([0,\infty);H^1_0(\Omega)\big)\cap W^{1,\frac{m+1}m}_\loc\big([0,\infty);X^\star+L^\frac2m(\Omega)\big),	\medskip \\
\big(g_\eps^m(u_\eps)\big)_{\eps>0} \text{ is bounded in } L^\infty_\loc\big([0,\infty);L^\frac2m(\Omega)\big),
\end{cases}
\end{gather}
and
\begin{gather}
\label{AlemthmsaH16}
\sup_{\eps>0}\vint_0^T\!\!\!\vint_\Omega\frac{|u_\eps(t,x)|^2}{(|u_\eps(t,x)|^2+\eps)^\frac{1-m}2}\d x\d t\le C(T),
\end{gather}
for any $T>0.$ The full correct version may be found in~\cite{arXiv}.
\end{rmk}

\begin{vproof}{of Proposition~\ref{propdep}.}
The embedding in \eqref{propdephyp} comes from Bégout and D\'iaz~\cite[Lemma~A.4]{MR4053613} $(1<p<\infty)$ and Theorem~\ref{thmembL2} $(p=1$ or $p=\infty).$ We make the difference between the two equations satisfied by $u$ and $\wt u,$ respectively. It follows from \eqref{lemVL2} that $u-\wt u$ satisfies the equation obtained in $L^1_\loc(0,\infty;X^\star).$ We take the $X^\star-X$ duality product with $\vi(u-\wt u).$ By Lemma~\ref{lemAmon}, Theorem~\ref{thmembL2}, Bégout and D\'iaz~\cite[Lemma~A.5]{MR4053613}, and Cauchy-Schwarz' inequality, we get,
\begin{gather*}
\frac12\frac\d{\d t}\|u-\wt u\|_{L^2(\Omega)}^2\le\|f-\wt f\|_{L^2(\Omega)}\|u-\wt u\|_{L^2(\Omega)},
\end{gather*}
almost everywhere on $(0,\infty).$ Integrating over $(s,t),$ we obtain \eqref{estthmweak}. Finally, we note that the strong solutions satisfy \eqref{propdephyp} with $p=1,$ and that \eqref{estthmweak} is stable by passing to the limit in $C\big([0,T];L^2(\Omega)\big)\times L^1(0,T;L^2(\Omega)),$ for any $T>0.$ By using \eqref{cv}, we then deduce that \eqref{estthmweak} still holds true for the weak solutions.
\medskip
\end{vproof}

\begin{vproof}{of Theorem~\ref{thmstrongaH2}.}
Let $f$ and $u_0$ be as in the theorem. We recall that we identify the Hilbert space $L^2(\Omega)$ with its own dual. It follows that the duality mapping is nothing else but the identity. As a consequence, and from Corollary~\ref{corAmon} and  Lemma~\ref{lemAmax}, the operator $(A,D(A))$ is $m$-accretive (i.e., maximal monotone) on $L^2(\Omega).$ Then, by Vrabie~\cite[Theorem~1.7.1]{MR1375237}), there exist a unique $u\in W^{1,\infty}_\loc\big([0,\infty);L^2(\Omega)\big)$ and $U\in L^\infty_\loc\big([0,\infty);H^{-2}(\Omega)\big)$ which satisfy \eqref{nls} in $\Dr^\p\big((0,\infty)\times\Omega\big),$ \eqref{strongaH23}, $u(0)=u_0,$ $-\vi\Delta u-\vi Vu+\mu\,U\in L^\infty_\loc(0,\infty;L^2(\Omega)).$ In addition, for almost every $t>0,$ $u(t)\in H^1_0(\Omega)\cap L^1(\Omega),$ $U(t)\in\ovl{B}_{L^\infty}(0,1),$ $(-\vi\Delta u+\mu\,U)(t)\in L^2(\Omega),$ and $U(t,x)=\frac{u(t,x)}{|u(t,x)|},$ for almost every $(t,x)\in(0,\infty)\times\Omega$ where $u(t,x)\neq0.$ Then \eqref{strongaH21} is an immediate consequence of \eqref{strongaH23}. Taking the $L^2$-scalar product of \eqref{nls} with $\vi u,$ we deduce from \cite[Lemma~A.5]{MR4053613}, and \cite[Lemma~4.4]{MR4503241} that Property~\ref{thmstrongaH23} holds true. Then, $u\in L^\infty_\loc\big([0,\infty);L^1(\Omega)\big),$ and applying the Cauchy-Schwarz inequality to \eqref{L2}, and integrating, we get \eqref{strongaH22}. Now, let us establish \eqref{strongaH24} and prove that $u\in C_\w\big([0,\infty);H^1_0(\Omega)\big).$ We take again the $L^2$-scalar product of \eqref{nls} with $-u$ and apply \cite[Lemma~4.4]{MR4503241}. By the Cauchy-Schwarz inequality, we infer
\begin{gather*}
\|\nabla u\|_{L^2(\Omega)}^2\le\left(\|u_t\|_{L^2(\Omega)}+\|Vu\|_{L^2(\Omega)}+\|f\|_{L^2(\Omega)}\right)\|u\|_{L^2(\Omega)},
\end{gather*}
and with help of (4.2) and (4.3) in \cite{MR4340780}, we have that for some $\gamma\in[0,1),$
\begin{gather*}
\|Vu\|_{L^2(\Omega)}\le C(N,\beta)\left(\|V_1\|_{L^\infty(\Omega)}+\|V_2\|_{L^{p_V}(\Omega)}^{2-\gamma}\right)\|u\|_{L^2(\Omega)}
		+\frac1{2\|u\|_{L^2(\Omega)}}\|\nabla u\|_{L^2(\Omega)}^2.
\end{gather*}
Putting together the two above estimates, we deduce that
\begin{gather}
\label{demthmstrongaH21}
\|\nabla u\|_{L^2(\Omega)}^2\le C\left(\|u_t\|_{L^2(\Omega)}+\left(\|V_1\|_{L^\infty(\Omega)}+\|V_2\|_{L^{p_V}(\Omega)}^{2-\gamma}\right)
  			\|u\|_{L^2(\Omega)}+\|f\|_{L^2(\Omega)}\right)\|u\|_{L^2(\Omega)},
\end{gather}
almost everywhere on $(0,\infty),$ where $C=C(N,\beta).$ In particular, \eqref{strongaH22}, \eqref{strongaH23}, \eqref{L2} and \eqref{demthmstrongaH21} give \eqref{strongaH24}. Now, since $u\in C\big([0,\infty);L^2(\Omega)\big),$ it follows from \eqref{strongaH22}, \eqref{strongaH24} and \eqref{weakcon} that $u\in C_\w\big([0,\infty);H^1_0(\Omega)\big).$ In particular, by \eqref{lemVL2}, we have that $Vu\in L^\infty_\loc\big([0,\infty);L^2(\Omega)\big).$ We claim that $U\in L^\infty\big((0,\infty)\times\Omega\big)$ and that $\|U\|_{L^\infty((0,\infty)\times\Omega)}\le1.$ We first have to show that $(t,x)\longmapsto U(t,x)$ is measurable. By the previous estimates, we only have that,
\begin{gather}
\label{demthmstrongaH22}
-\vi\Delta u+\mu\,U\in L^\infty_\loc\big([0,\infty);L^2(\Omega)\big),
\end{gather}
but we do not have any information on $\Delta u$ and $U,$ separately. Let $\Omega^\p\subset\Omega$ be any bounded open subset. Since for almost every $t>0,$ $U(t)\in\ovl{B}_{L^\infty}(0,1),$ it follows that for any $T>0,$ there exist $C>0$ and $N_0\subset(0,T)$ such that $|N_0|=0,$ and for any $t\in(0,T)\setminus N_0,$ $\|\Delta u(t)_{|\Omega^\p}\|_{L^2(\Omega^\p)}\le C.$ But, $\Delta u(\,.\,)_{|\Omega^\p}\in C_\w\big([0,\infty);H^{-1}(\Omega^\p)\big)$ and so by \eqref{weakcon}, $\Delta u(\,.\,)_{|\Omega^\p}\in C_\w\big([0,\infty);L^2(\Omega^\p)\big).$ It follows from \eqref{demthmstrongaH22} that
\begin{gather}
\label{demthmstrongaH23}
U\in  L^\infty_\loc\big([0,\infty);L^2_\loc(\Omega)\big)\inj L^2_\loc\big((0,\infty)\times\Omega\big).
\end{gather}
In particular, $U:(0,\infty)\times\Omega\tends\C$ is measurable and for almost every $(t,x)\in(0,\infty)\times\Omega,$ $|U(t,x)|\le\|U(t)\|_{L^\infty(\Omega)}\le1.$ Hence, $\|U\|_{L^\infty((0,\infty)\times\Omega)}\le1.$ Furthermore, by \eqref{demthmstrongaH23}, we infer that \eqref{nls} makes sense in $L^\infty_\loc\big([0,\infty);L^2_\loc(\Omega)\big)\cap L^\infty_\loc\big([0,\infty);H^{-1}(\Omega)\big).$ Finally, Property~\ref{thmstrongaH24} follows easily from Properties~\ref{thmstrongaH21}--\ref{thmstrongaH23} and \eqref{nls}. The theorem is proved.
\medskip
\end{vproof}

\begin{vproof}{of Theorem~\ref{thmweak}.}
Let $f\in L^1_\loc([0,\infty);L^2(\Omega))$ and $u_0\in L^2(\Omega).$ Let $(\vphi_n,f_n)_{n\in\N}\subset\Dr(\Omega)\times W^{1,1}_\loc([0,\infty);L^2(\Omega))$ be such that,
\begin{gather}
\label{demthmweak1}
\forall T>0, \; (\vphi_n,f_n)\xrightarrow[n\tends\infty]{L^2(\Omega)\times L^1(0,T;L^2(\Omega))}(u_0,f).
\end{gather}
Finally, for each $n\in\N,$ let $(u_n,U_n)$ be the unique $H^1_0$-solution to \eqref{nls}--\eqref{nlsb} with $u(0)=\vphi_n$ given by Theorem~\ref{thmstrongaH2}. It follows from \eqref{estthmweak} that for any $T>0,$ $(u_n)_{n\in\N}$ is a Cauchy sequence of the complete space $C\big([0,T];L^2(\Omega)\big).$ It follows that there exists $u\in C\big([0,\infty);L^2(\Omega)\big)$ such that,
\begin{gather}
\label{demthmweak2}
\forall T>0, \; u_n\xrightarrow[n\tends\infty]{C([0,T];L^2(\Omega))}u.
\end{gather}
Since $(U_n)_{n\in\N}$ is bounded in $L^\infty\big((0,\infty)\times\Omega)$ by one, there exists $U\in L^\infty\big((0,\infty)\times\Omega)$ satisfying \eqref{defsol22} such that, up to a subsequence, \eqref{U} holds true. Then $u$ is a weak solution to \eqref{nls}--\eqref{u0}, and by Proposition~\ref{propdep}, this solution is unique. Let $t\ge s\ge0.$ Using Theorem~\ref{thmstrongaH2}, we obtain that for any $n\in\N,$
\begin{gather}
\label{demthmweak3}
\frac12\|u_n(t)\|_{L^2(\Omega)}^2+\mu\dsp\vint_s^t\|u_n(\sigma)\|_{L^1(\Omega)}\d\sigma
=\frac12\|u_n(s)\|_{L^2(\Omega)}^2+\Im\dsp\iint\limits_{s\;\Omega}^{\text{}\;\;t}f_n(\sigma,x)\,\ovl{u_n(\sigma,x)}\,\d x\,\d\sigma,
\end{gather}
If $|\Omega|<\infty$ then $L^2(\Omega)\inj L^1(\Omega).$ We then use \eqref{demthmweak1}--\eqref{demthmweak2} to pass to the limit in \eqref{demthmweak3} from which \eqref{L1}--\eqref{L2+} follows. If $|\Omega|=\infty$ then by \eqref{demthmweak2}, there exists a subsequence that we still denote by $(u_n)_{n\in\N}$ such that,
\begin{gather}
\label{demthmweak4}
u_n\xrightarrow[n\tends\infty]{\text{a.e.\,in\:}(0,\infty)\times\Omega}u.
\end{gather}
Using \eqref{demthmweak1}, \eqref{demthmweak2} and \eqref{demthmweak4}, we deduce \eqref{L1}--\eqref{L2+} from \eqref{demthmweak3} and Fatou's Lemma.
\medskip
\end{vproof}

\section{Proofs of the finite time extinction and the asymptotic behavior theorems}
\label{proofext}

\begin{vproof}{of Theorems~\ref{thmextN1}.}
We first consider the assumption \eqref{f}. By \eqref{L2} and Hölder's inequality and \eqref{f}, we have for almost every $t\in(T_0,\infty),$
\begin{gather}
\label{edo}
\frac12\frac\d{\d t}\|u(t)\|_{L^2(\Omega)}^2+\lambda\|u(t)\|_{L^1(\Omega)}\le0,
\end{gather}
where $\lambda\eqdef\mu-\|f\|_{L^\infty((T_0,\infty)\times\Omega)}>0.$ By the Gagliardo-Nirenberg inequality, there exists $C_\GN=C(N)$ such that for almost every $t>0,$
\begin{gather}
\label{GN}
\|u(t)\|_{L^2(\Omega)}^\frac{N+2}2\le C_\GN\|u(t)\|_{L^1(\Omega)}\|\nabla u(t)\|_{L^2(\Omega)}^\frac{N}2.
\end{gather}
By Theorem~\ref{thmstrongaH2}, $u\in L^\infty(0,\infty;H^1_0(\Omega)).$ Then, setting for any $t\ge0,$ $y(t)=\|u(t)\|_{L^2(\Omega)}^2,$ and using \eqref{GN} in \eqref{edo}, we get that for almost every $t>0,$
\begin{gather*}
y^\p(t)+C\|\nabla u\|_{L^\infty(0,\infty;L^2(\Omega))}^{-\frac{N}2}y(t)^\frac{N+2}4\le0.
\end{gather*}
Integrating over $(0,T)$ the above estimate, we obtain \eqref{01}--\eqref{thmrtdH12} (see also \cite[Lemma~5.1]{MR4053613}). It remains to show the last property on the instantaneous extinction time. By \eqref{L2}, \eqref{GN}, and Young's inequality, we have for almost every $t>0,$
\begin{gather}
\label{demthmextN1}
y^\p(t)+2\,\mu\,C_\GN^{-1}\,\|\nabla u\|_{L^\infty(0,\infty;L^2(\Omega))}^{-\frac12}y(t)^\frac34\le3\|f(t)\|_{L^2(\Omega)}^3+y(t)^\frac34.
\end{gather}
Let $\delta=\dfrac34,$ and
\begin{gather*}
x_\star=(\delta(1-\delta)T_0)^\frac1{1-\delta}, \quad y_\star=(\delta^\delta(1-\delta))^\frac1{1-\delta},	\\
\eps_1=\min\left\{(\delta(1-\delta))^\frac1{2(1-\delta)},\left(\frac{y_\star}3\right)^\frac13\right\}.
\end{gather*}
Let $M>0.$ Assume that $\|f\|_{W^{1,1}(0,\infty;L^2(\Omega))}\le M.$ By Sobolev' embedding, this entails that $\|f(0)\|_{L^2(\Omega)}\le C(M).$ Let $\omega=\{x\in\Omega; u_0(x)\neq0\}.$ We first need a refinement of $B(t)$ in \eqref{strongaH23}. Actually, we have for any $t>0$ by a result in Barbu~\cite{MR2582280} (Theorem~4.4, p.141),
\begin{gather*}
\left\|u_t\right\|_{L^\infty(0,t;L^2(\Omega))}\le
\wt B(t)\eqdef|\Delta u_0+Vu_0+\vi\mu U_0-f(0)|_{L^2(\Omega)}+\int_0^t\|f^\p(\sigma)\|_{L^2(\Omega)}\d\sigma,
\end{gather*}
where $|\Delta u_0+Vu_0+\vi\mu U_0-f(0)|_{L^2(\Omega)}=\min\big\{\|v\|_{L^2(\Omega)};v\in\vi Au_0-f(0)\big\}.$ (The multi-valued operator $(A,D(A))$ is defined at the Section~\ref{proofexi}.) Furthermore, this minimum is unique (see the lines before Barbu~\cite[Proposition~3.5, p.101]{MR2582280}). It follows that $B(t)$ may be replaced with $\wt B(t)$ in the whole Theorem~\ref{thmstrongaH2}. Let us define $V_0$ in $\Omega$ by $V_0=|u_0|^{-1}u_0,$ in $\omega,$ and $V_0=0,$ in $\omega^\co.$ Since $u_0\in H^2_\loc(\Omega),$ we have in particular that $\Delta u_0=0,$ almost everywhere in $\omega^\co.$ Then,
\begin{align*}
\sup_{t>0}\wt B(t)\le	&	\; \|\Delta u_0+Vu_0+\vi\mu V_0-f(0)\|_{L^2(\Omega)}+\int_0^\infty\|f^\p(\sigma)\|_{L^2(\Omega)}\d\sigma, 		\\
			  \le	&	\; \|\Delta u_0+Vu_0+\vi\mu|u_0|^{-1}u_0\|_{L^2(\omega)}+\|f(0)\|_{L^2(\Omega)}
												+\int_0^\infty\|f^\p(\sigma)\|_{L^2(\Omega)}\d\sigma.
\end{align*}
As a consequence, if $\|u_0\|_{L^2(\Omega)}\le\eps_1$ and $\|\nabla u_0\|_{L^2(\Omega)}+\left\|\Delta u_0+\vi\mu\frac{u_0}{|u_0|}\right\|_{L^2(\omega)}\le M,$ then it follows from \eqref{lemVL2} and the above that
\begin{align*}
\sup_{t>0}\wt B(t)\le	&	\; \left\|\Delta u_0+\vi\mu\frac{u_0}{|u_0|}\right\|_{L^2(\omega)}+\|Vu_0\|_{L^2(\Omega)}+C(M)	\\
			  \le	&	\; C(\|V\|_{L^{p_V}+L^\infty},M).
\end{align*}
By \eqref{strongaH24}, there exists a positive constant $\eps_\star=\eps_\star(\|V\|_{L^{p_V}+L^\infty},M,\mu)\le\eps_1$ such that if $\|u_0\|_{L^2(\Omega)}+\|f\|_{L^1(0,\infty;L^2(\Omega))}\le\eps_\star$ then $\mu\,C_\GN^{-1}\,\|\nabla u\|_{L^\infty((0,\infty);L^2(\Omega))}^{-\frac12}\ge1.$ Now, assume that $\|f(t)\|_{L^2(\Omega)}\le\eps_\star\big(T_0-t\big)_+,$ for any $t\ge0.$ Then \eqref{demthmextN1} becomes,
\begin{gather*}
\text{for almost every } t>0, \;\; y^\p(t)+y(t)^\frac34\le y_\star\big(T_0-t\big)_+^\frac\delta{1-\delta}.
\end{gather*}
Finally, if $\|u_0\|_{L^2(\Omega)}\le\eps_\star T_0^2$ then $y(0)\le x_\star,$ and an appeal to \cite[Lemma~5.2]{MR4053613} gives that for any $t\ge T_0,$ $y(t)=0.$ The theorem is proved.
\medskip
\end{vproof}

\begin{vproof}{of Theorem~\ref{thm0w}.}
By density and \eqref{estthmweak}, it is sufficient to consider the case in which $u_0\in\Dr(\Omega)$ and $f\in\Dr\big([0,\infty);L^2(\Omega)\big).$ Then Theorem~\ref{thmextN1} yields the desired result.
\medskip
\end{vproof}

\section{Solutions with lower regularity}
\label{secaH1}

Due to problems of measurability of functions with values in non-separable spaces (Section~\ref{secnm}), we are unable to build solutions to \eqref{nls} under the mere ``natural condition" $u_0\in H^1_0(\Omega),$ and $f$ which satisfies \eqref{secnmf}. The closest result is the following.

\begin{thm}
\label{thmstrongaH1}
Let Assumption~$\ref{ass1}$ be fulfilled. Let $(Y_n)_{n\in\N_0}$ be any $L^1$-approximating sequence of RNP-spaces $($Definition~$\ref{defRNP}).$ Assume that $V\in W^{1,\infty}(\Omega;\R)+W^{1,p_V}(\Omega;\R),$ where $p_V$ is given by~\eqref{pV} and let
\begin{gather}
\label{thmstrongaH11}
f\in L^1_\loc\big([0,\infty);H^1_0(\Omega)\big)\cap L^p_\loc\big([0,\infty);H^{-1}(\Omega)+L^\infty(\Omega)\big),
\end{gather}
for some $1<p\le\infty.$ Then for any $u_0\in H^1_0(\Omega),$ there exist $u$ and $U$ satisfying~\eqref{defsol22}--\eqref{defsol23} such that for any $n\in\N_0,$
\begin{gather}
\begin{cases}
\label{thmstrongaH12}
u\in C_\w\big([0,\infty);H^1_0(\Omega)\big)\cap W^{1,p}_\loc\big([0,\infty);H^{-1}(\Omega)+Y_n^\star\big),	\medskip \\
u\in L^1_\loc\big([0,\infty);L^1(\Omega)\big),
\end{cases}
\end{gather}
and such that the pair $(u,U)$ is a solution to~\eqref{defsol21} in $L^p_\loc(0,\infty;H^{-1}(\Omega)+Y_n^\star)\inj\Dr^\p\big((0,\infty)\times\Omega\big).$ Furthermore, $u$ verifies $u(0)=u_0,$ and
\begin{gather}
\label{thmstrongaH13}
\|u(t)\|_{H^1_0(\Omega)}\le\left(\|u_0\|_{H^1_0(\Omega)}+\vint_0^t\|f(s)\|_{H^1_0(\Omega)}\d s\right)e^{C\|\nabla V\|_{L^\infty+L^{p_V}}t},
\end{gather}
for any $t\ge0,$ where $C=C(N,\beta).$ In addition, if $\nabla V=0$ then we have that,
\begin{gather}
\label{thmstrongaH14}
\|\nabla u(t)\|_{L^2(\Omega)}\le\|\nabla u_0\|_{L^2(\Omega)}+\vint_0^t\|\nabla f(s)\|_{L^2(\Omega)}\d s,
\end{gather}
for any $t\ge0.$ Moreover, there exists $N_0\subset(0,\infty)$ with $|N_0|=0$ such that,
\begin{gather}
\label{thmstrongaH15}
\forall t\in(0,\infty)\setminus N_0, \; u^\p(t)\in H^{-1}(\Omega)+L^\infty(\Omega).
\end{gather}
In particular, $u$ solves~\eqref{defsol21} in $ H^{-1}(\Omega)+L^\infty(\Omega),$ for almost every $t>0.$ Finally, if $p=\infty$ then for any $T>0,$ there exists $C(T)>0$ such that,
\begin{gather}
\label{thmstrongaH16}
\|u^\p(t)\|_{H^{-1}(\Omega)+L^\infty(\Omega)}\le C(T),
\end{gather}
for any $t\in(0,T)\cap N_0^\co.$
\end{thm}

\begin{rmk}
\label{rmkthmstrongaH1}
Below are some comments about Theorem~\ref{thmstrongaH1}.
\begin{enumerate}
\item
\label{rmkthmstrongaH11}
We do not know whether $u\in C\big([0,\infty);L^2(\Omega)\big),$ and we do not know whether $u^\p:(0,\infty)\tends H^{-1}(\Omega)+L^\infty(\Omega)$ is measurable.
\item
\label{rmkthmstrongaH12}
Let $f$ satisfy \eqref{thmstrongaH11} with $p=\infty,$ and let $X=H^1_0(\Omega)\cap L^1(\Omega).$ If $u^\p:(0,\infty)\tends X^\star$ is measurable then it follows from \eqref{thmstrongaH12}, \eqref{thmstrongaH15} and \eqref{thmstrongaH16} that
\begin{gather*}
u\in L^1_\loc\big([0,\infty);X\big)\cap W^{1,\infty}_\loc\big([0,\infty);X^\star\big),
\end{gather*}
so that $u$ becomes an $H^1_0$-solution.
\item
\label{rmkthmstrongaH13}
Let $f$ satisfy \eqref{thmstrongaH11} with $p=\infty.$ Actually, we may get a little bit more information from~\eqref{thmstrongaH15}--\eqref{thmstrongaH16}. Let $n\in\N_0,$ $X_n=H^1_0(\Omega)\cap Y_n,$ and $X=H^1_0(\Omega)\cap L^1(\Omega).$ By~\eqref{thmstrongaH12}, we have that
\begin{gather}
\label{diffYn*}
\lim_{h\to0}\left\|\frac{u(t+h)-u(t)}h-u^\p(t)\right\|_{X_n^\star}=0,
\end{gather}
for almost every $t>0.$ But $X_n\inj X,$ with dense embedding, and $X$ is separable. It then follows from \eqref{thmstrongaH16} and \eqref{diffYn*} that,
\begin{gather}
\label{ufmX*}
\frac{u(t+h)-u(t)}h\underset{h\to0}{-\!\!\!-\!\!\!-\!\!\!-\!\!\!\weak}u^\p(t), \text{ in the weak}\!\star \text{topology } \sigma(X^\star,X),
\end{gather}
for almost every $t>0,$ so that $u^\p:(0,\infty)\tends X^\star$ is weak$\star$-measurable. Thus, with \eqref{thmstrongaH16}, we conclude that $u^\p$ is Gel'fand integrable: for any finite-measure set $E\subset(0,\infty)$ and $\vphi\in X,$
\begin{gather*}
\left\langle G-\int_E u^\p(t)\d t,\vphi\right\rangle_{X^\star,X}=\int_E\langle u^\p(t),\vphi\rangle_{X^\star,X}\d t,
\end{gather*}
where $G-\dsp\int_E u^\p(t)\d t$ denotes the Gel'fand integral over $E$ (Gel$^\p$fand~\cite{gelfand}). In addition, there exists a sequence $(s_n)_{n\in\N}$ of simple functions such that for any $\vphi\in X,$
\begin{gather*}
\lim_{n\to\infty}\langle s_n(t),\vphi\rangle_{X^\star,X}=\langle u^\p(t),\vphi\rangle_{X^\star,X},
\end{gather*}
for almost every $t\in E.$ Note that the null set $N_\vphi$ on which the convergence fails may vary with $\vphi$ (while if $u^\p:(0,\infty)\tends X^\star$ was measurable, $N_\vphi\equiv N$ would not depend on $\vphi,$ and $u^\p:(0,\infty)\tends X^\star$ would be locally Bochner integrable). By Pettis' Theorem, $u^\p:(0,\infty)\tends X^\star$ is measurable if, and only if, there exists a null set $N_0\subset(0,\infty)$ for which $u^\p((0,\infty)\setminus N_0)$ is separable in $X^\star.$ For more details about weak$\star$-measurable functions, the Gel'fand integral and its properties, see Gel$^\p$fand~\cite{gelfand}, Hashimoto~\cite{MR707846}, Hashimoto and Oharu~\cite{MR707185}, and Schwartz~\cite{MR0448075}.
\item
\label{rmkthmstrongaH14}
As noted just above, $u^\p:(0,\infty)\tends X^\star$ is weak$\star$-measurable, where $X=H^1_0(\Omega)\cap L^1(\Omega).$ It seems very hard to answer to the question of the measurability of $u^\p.$ Indeed, the first way, is to know whether $u^\p$ is essentially separably valued in $X^\star$ (Pettis' Theorem). The second method would be to show that $u^\p$ is weakly$\star$ equivalent to a measurable function, namely, to prove that there exists a measurable function $v:(0,\infty)\tends X^\star$ such that for any $\vphi\in X,$
\begin{gather*}
\langle u^\p(t),\vphi\rangle_{X^\star,X}= \langle v(t),\vphi\rangle_{X^\star,X},
\end{gather*}
for almost every $t>0.$ A priori, the null set on which the equality fails may vary with $\vphi,$ but actually, it does not with help of the separability of $X.$ Therefore, we would have $u^\p(t)=v(t),$ for almost every $t>0,$ so that $u^\p$ would be measurable. But showing that $u^\p$ is weakly$\star$ equivalent to a measurable function requires to be able to identify weakly compact sets of $X^\star,$ namely in the weak topology $\sigma(X^\star,X^{\star\star})$ (Uhl~\cite{MR466473}, Edgar~\cite{MR487448}). A third method would consist to see the Gel'fand integral of $u^\p$ as a vector measure and to look for make it represent by a measurable function with help of the Radon-Nikod\'ym Theorem (Rieffel~\cite{MR222245}, Maynard~\cite{MR385521}, Diestel and Uhl~\cite{MR399852}, Schwartz~\cite{MR0500143}). It is interesting to note that all these theorems furnish a necessary and sufficient condition to measurability. Unfortunately, they are too general and too abstract to be applied in our case (such as the $\sigma$-dentability, for the third method), and the only handy tools about measurability of the vector measure theory are reflexivity and separability.
\end{enumerate}
\end{rmk}

\begin{lem}
\label{lemthmsaH1}
We use the notations of Theorem~$\ref{thmstrongaH1}.$ Let the hypotheses of Theorem~$\ref{thmstrongaH1}$ be fulfilled with $p<\infty.$ Let $(f_\eps)_{\eps>0}\subset\Dr\big([0,\infty);H^1_0(\Omega)\big)$ and $(\vphi_\eps)_{\eps>0}\subset\Dr(\Omega)$ be such that,
\begin{gather}
\label{lemthmsaH11}
\begin{cases}
\vphi_\eps\xrightarrow[\eps\searrow0]{H^1_0(\Omega)}u_0,	\\
f_\eps\xrightarrow[\eps\searrow0]{L^1(0,T;H^1_0)\cap  L^p(0,T;H^{-1}+L^\infty)}f.
\end{cases}
\end{gather}
for any $T>0.$ Then for any $\eps>0,$ there exists a unique solution
\begin{gather}
\label{lemthmsaH12}
u_\eps\in C_\w\big([0,\infty);H^1_0(\Omega)\big)\cap W^{1,\infty}_\loc\big([0,\infty);L^2(\Omega)\big),
\end{gather}
to
\begin{gather}
\label{nlse}
\vi\frac{\partial u_\eps}{\partial t}+\Delta u_\eps+V(x)u_\eps+\vi\mu g_\eps(u_\eps)=f_\eps(t,x),	\text{ in } L^2(\Omega),
\end{gather}
for almost every $t>0,$ such that $u_\eps(0)=\vphi_\eps.$ Furthermore, for any $\eps>0,$
\begin{gather}
\label{lemthmsaH13}
\|u_\eps(t)\|_{H^1_0(\Omega)}
\le\left(\|\vphi_\eps\|_{H^1_0(\Omega)}+\vint_0^t\|f_\eps(s)\|_{H^1_0(\Omega)}\d s\right)e^{C\|\nabla V\|_{L^\infty+L^{p_V}}t},
\end{gather}
for any $t\ge0,$ where $C=C(N,\beta),$ and if $\nabla V=0$ then,
\begin{gather}
\label{lemthmsaH14}
\|\nabla u_\eps(t)\|_{L^2(\Omega)}\le\|\nabla\vphi_\eps\|_{L^2(\Omega)}+\vint_0^t\|\nabla f_\eps(s)\|_{L^2(\Omega)}\d s,
\end{gather}
for any $t\ge0.$ Finally,
\begin{gather}
\label{lemthmsaH15}
\begin{cases}
(u_\eps)_{\eps>0} \text{ is bounded in }
L^\infty_\loc\big([0,\infty);H^1_0(\Omega)\big)\cap W^{1,p}_\loc\big([0,\infty);X_n^\star\big),	\medskip \\
\big(g_\eps(u_\eps)\big)_{\eps>0} \text{ is bounded in } L^\infty(0,\infty;X_n^\star),
\end{cases}
\end{gather}
for any $n\in\N_0,$ where $X_n=H^1_0(\Omega)\cap Y_n,$ and
\begin{gather}
\label{lemthmsaH16}
\sup_{\eps>0}\vint_0^T\!\!\!\vint_\Omega\frac{|u_\eps(t,x)|^2}{(|u_\eps(t,x)|^2+\eps)^\frac12}\d x\d t\le C(T),
\end{gather}
for any $T>0.$
\end{lem}

\begin{proof*}
Let the assumptions of the Lemma be fulfilled. By \cite[Corollary~5.11]{MR4340780} and Barbu~\cite[Theorem~4.5, p.141]{MR2582280}, there exists a unique solution $u_\eps\in W^{1,\infty}_\loc\big([0,\infty);L^2(\Omega)\big)$ to \eqref{nlse} such that $u_\eps(t)\in H^1_0(\Omega)$ and $\Delta u_\eps(t)\in L^2(\Omega),$ for almost every $t>0.$ Moreover, $u_\eps(0)=\vphi_\eps.$ Taking the $L^2$-scalar product of \eqref{nlse} with $-u_\eps,$ applying Cauchy-Schwarz' inequality, and then (4.2) and (4.3) of \cite{MR4340780}, we get that $u_\eps\in L^\infty_\loc\big([0,\infty);H^1_0(\Omega)\big).$ Therefore, by \eqref{weakcon} we get \eqref{lemthmsaH12}. By \eqref{lemVL2} and \eqref{nlse}, it follows that $\Delta u_\eps\in L^\infty_\loc\big([0,\infty);H^1_0(\Omega)\big).$ Taking the $L^2$-scalar product of \eqref{nlse} with $-\vi\Delta u_\eps,$ it then follows from \cite{MR4053613} ((6.8) and Lemma~A.5), \eqref{lemVL2} and a density argument that for almost every $s>0,$
\begin{gather}
\label{demlemthmsaH11}
\frac12\frac{\d}{\d t}\|\nabla u_\eps(s)\|_{L^2(\Omega)}^2\le\|\nabla f_\eps(s)\|_{L^2(\Omega)}\|\nabla u_\eps(s)\|_{L^2(\Omega)}
								+C\|\nabla V\|_{L^\infty(\Omega)+L^{p_V}(\Omega)}\|u_\eps(s)\|_{H^1_0(\Omega)}^2.
\end{gather}
Therefore, \eqref{lemthmsaH14} follows by integrating \eqref{demlemthmsaH11}. Taking again the $L^2$-scalar product of \eqref{nlse} with $\vi u_\eps,$ we get that,
\begin{gather}
\label{demlemthmsaH12}
\frac12\frac{\d}{\d t}\|u_\eps(s)\|_{L^2(\Omega)}^2+\mu\vint_\Omega\frac{|u_\eps(s,x)|^2}{(|u_\eps(s,x)|^2+\eps)^\frac12}\d x
\le\|f_\eps(s)\|_{L^2(\Omega)}\|u_\eps(s)\|_{L^2(\Omega)},
\end{gather}
for almost every $s>0.$ Summing \eqref{demlemthmsaH11} with \eqref{demlemthmsaH12}, we get for almost every $s>0,$
\begin{gather*}
\frac12\frac{\d}{\d t}\|u_\eps(s)\|_{H^1_0(\Omega)}^2\le\|f_\eps(s)\|_{H^1_0(\Omega)}\|u_\eps(s)\|_{H^1_0(\Omega)}
+C\|\nabla V\|_{L^\infty(\Omega)+L^{p_V}(\Omega)}\|u_\eps(s)\|_{H^1_0(\Omega)}^2.
\end{gather*}
Integrating the above and applying Gronwall's Lemma, we get \eqref{lemthmsaH13}, implying that
\begin{gather}
\label{demlemthmsaH13}
(u_\eps)_{\eps>0} \text{ is bounded in } L^\infty_\loc\big([0,\infty);H^1_0(\Omega)\big).
\end{gather}
So, integrating \eqref{demlemthmsaH12}, we get \eqref{lemthmsaH16}. Let $n\in\N_0.$ Since $u^\p_\eps\in L^\infty_\loc\big([0,\infty);L^2(\Omega)\big),$ it follows from \eqref{nlse} that both $u^\p_\eps$ and $g_\eps(u_\eps)$ are measurable $[0,\infty)\tends X_n^\star.$ Therefore, since $|g_\eps(u_\eps)|\le1,$ almost everywhere in $(0,\infty)\times\Omega,$ we get \eqref{lemthmsaH15} by \eqref{demlemthmsaH13}, \eqref{lemVL2}, Lemma~\ref{lemY*}, Corollary~\ref{corZY}, and again \eqref{nlse}. This ends the proof of the lemma.
\medskip
\end{proof*}

\begin{lem}
\label{lemthmstrongaH23}
We use the notations of Lemma~$\ref{lemthmsaH1}.$ Under the hypotheses of Lemma~$\ref{lemthmsaH1},$ there exist $u$ and $U\in L^\infty(0,\infty;X_n^\star)$ satisfying~\eqref{defsol22}--\eqref{defsol23} such that,
\begin{gather}
\label{lemthmstrongaH23a}
u\in C_\w\big([0,\infty);H^1_0(\Omega)\big)\cap W^{1,p}_\loc\big([0,\infty);X_n^\star\big),	\\
\label{lemthmstrongaH23b}
u\in  L^1_\loc\big([0,\infty);L^1(\Omega)\big),
\end{gather}
for any $n\in\N_0,$ and a positive sequence $\eps_\ell\searrow0,$ as $\ell\tends\infty,$ such that
\begin{align}
\label{lemthmstrongaH23c}
&	u_{\eps_\ell}(t)\underset{\ell\to\infty}{-\!\!\!-\!\!\!-\!\!\!-\!\!\!\weak}u(t) \text{ in } H^1_0(\Omega)_\w, \; \forall t\ge0,		\\
\label{lemthmstrongaH23d}
&	u_{\eps_\ell}\xrightarrow[n\to\infty]{\text{a.e.\,in }(0,\infty)\times\Omega}u,									\\
\label{lemthmstrongaH23e}
&	g_{\eps_\ell}\big(u_{\eps_\ell}\big)\underset{\ell\to\infty}{-\!\!\!-\!\!\!-\!\!\!-\!\!\!\weak}U, \text{ in } L^\infty(0,\infty;X_n^\star)_{\w\star}.
\end{align}
\end{lem}

\begin{proof*}
Let $n\in\N.$ By \eqref{lemthmsaH12}, \eqref{lemthmsaH15}, Cazenave~\cite{MR2002047} (Proposition~1.1.2(i), p.2, and Remark~1.3.13(ii), p.12)) and the diagonal procedure, there exist a positive sequence $\eps_\ell\searrow0,$ as $\ell\tends\infty,$ and $u$ satisfying \eqref{lemthmstrongaH23a} and \eqref{lemthmstrongaH23c}. Let $T>0$ and $\Omega^\p\subset\Omega$ be a smooth bounded open subset of $\R^N.$ Let $(Y_n(\Omega^\p))_{n\in\N_0}$ be an $L^1(\Omega^\p)$-approximating sequence of RNP-spaces (Lemma~\ref{lemY}). By the Rellich-Kondrachov compactness Theorem,
\begin{gather}
\label{demlemthmsaH23a}
H^1(\Omega^\p)\underset{\text{compact}}{\inj} L^2(\Omega^\p)\inj H^{-1}(\Omega^\p)+Y_n(\Omega^\p)^\star.
\end{gather}
By \eqref{lemthmsaH12} and \eqref{nlse}, $u^\p_\eps$ and $g_\eps(u_\eps)$ are measurable $[0,\infty)\tends H^{-1}(\Omega^\p)+Y_n(\Omega^\p)^\star.$ Therefore, since $|g_\eps(u_\eps)|\le1,$ almost everywhere on $(0,\infty)\times\Omega,$ we get by \eqref{lemthmsaH15}, \eqref{lemVL2}, the embedding $L^\infty(\Omega^\p)\inj Y_n(\Omega^\p)^\star,$ and again \eqref{nlse}, that
\begin{gather}
\left\{
\begin{split}
\label{demlemthmsaH23b}
&	(u_\eps)_{\eps>0} \text{ is bounded in } L^\infty_\loc\big([0,\infty);H^1(\Omega^\p)\big)	\\
&	\; \text{and in } W^{1,p}_\loc\big([0,\infty);H^{-1}(\Omega^\p)+Y_n(\Omega^\p)^\star\big).
\end{split}
\right.
\end{gather}
By \eqref{lemthmstrongaH23c}, \eqref{demlemthmsaH23a}--\eqref{demlemthmsaH23b} and compactness (Simon~\cite[Corollary~4, p.85]{MR916688}), we obtain that
\begin{gather*}
u\in C\big([0,T];L^2(\Omega^\p)\big) \text{ and } \lim_{\ell\to\infty}\|u_{\eps_\ell}-u\|_{C([0,T];L^2(\Omega^\p))}=0.
\end{gather*}
Hence, $u_{\eps_\ell}\xrightarrow[\ell\to\infty]{L^2_\loc((0,\infty)\times\Omega)}u,$ since $T$ and $\Omega^\p$ are arbitrary. Up to a subsequence, we get \eqref{lemthmstrongaH23d}. We then obtain \eqref{lemthmstrongaH23b} from \eqref{lemthmsaH16}, \eqref{lemthmstrongaH23d} and Fatou's Lemma. Now, we note that by \eqref{injDp},
\begin{gather*}
L^\infty(0,\infty;X_n^\star)\inj\Dr^\p\big((0,\infty)\times\Omega\big).
\end{gather*}
In addition, by Corollary~\ref{corZY}, $X_n$ is separable and reflexive, so that $L^\infty(0,\infty;X_n^\star)$ is the dual space of the separable space $L^1(0,\infty;X_n)$ (Section~\ref{secfunana}). Finally, $\vsup_{\eps>0}\|g_\eps(u_\eps)\|_{L^\infty((0,\infty)\times\Omega)}\le1,$ so it follows from the above and \eqref{lemthmsaH15} that, up to a subsequence, there exists
\begin{gather*}
U\in L^\infty(0,\infty;X_n^\star)\cap L^\infty\big((0,\infty)\times\Omega\big),
\end{gather*}
satisfying~\eqref{defsol22}, \eqref{lemthmstrongaH23e}, and
\begin{gather*}
g_{\eps_\ell}\big(u_{\eps_\ell}\big)\underset{\ell\to\infty}{-\!\!\!-\!\!\!-\!\!\!-\!\!\!\weak}U, \text{ in } L^\infty\big((0,\infty)\times\Omega\big)_{\w\star}.
\end{gather*}
To conclude, it remains to show that $U$ satisfies \eqref{defsol23}. This may be done by repeating the proof of Lemma~\ref{lemAmax}, from \eqref{lemAmax5} to \eqref{lemAmax7}. The proof of the lemma is complete.
\medskip
\end{proof*}

\begin{vproof}{of Theorems~\ref{thmstrongaH1}.}
We use the notations of Lemma~$\ref{lemthmsaH1}.$ We first assume that $p<\infty.$ Let $u$ be given by Lemma~\ref{lemthmstrongaH23}. Let $n\in\N,$ let $\vphi\in X_n$ and $\psi\in C^1_\co\big((0,\infty);\R\big).$ Let $T>0$ be such that $\supp\psi\in(0,T).$ Let $(\eps_\ell)_{\ell\in\N}$ be given by Lemma~\ref{lemthmstrongaH23}. It follows from \eqref{nlse} and (4.5) in \cite{MR4340780}, that for any $\ell\in\N,$
\begin{gather*}
\vint_0^\infty\left\langle\vi\frac{\partial u_{\eps_\ell}}{\partial t}+\Delta u_{\eps_\ell}+Vu_{\eps_\ell}
	+\vi\mu g_{\eps_\ell}\big(u_{\eps_\ell}\big),\vphi\right\rangle_{X_n^\star,X_n}\psi(t)\,\d t=\vint_0^\infty\big\langle f_{\eps_\ell}(t),\vphi\big\rangle_{X_n^\star,X_n}\psi(t)\,\d t,
\end{gather*}
and so,
\begin{align*}
	&	\; \vint_0^T\Big(\left\langle-\vi u_{\eps_\ell},\vphi\right\rangle_{L^2(\Omega),L^2(\Omega)}\psi^\p(t)
			-\left\langle\nabla u_{\eps_\ell},\nabla\vphi\right\rangle_{L^2(\Omega),L^2(\Omega)}\psi(t)
			+\left\langle u_{\eps_\ell},V\vphi\right\rangle_{L^2(\Omega),L^2(\Omega)}\psi(t)	\Big)\d t							\\
   +	&	\; \left\langle\vi\mu g_{\eps_\ell}\big(u_{\eps_\ell}\big),\psi\vphi\right\rangle_{L^\infty((0,\infty);X_n^\star);L^1((0,\infty);X_n)}
			=\vint_0^T\big\langle f_{\eps_\ell}(t),\vphi\big\rangle_{L^2(\Omega),L^2(\Omega)}\psi(t)\,\d t.
\end{align*}
By \eqref{lemthmsaH11}, \eqref{lemthmstrongaH23c}, \eqref{lemthmstrongaH23e}, the dominated convergence Theorem, and Hölder's inequality, we can pass to the limit in the above equality to obtain $u(0)=u_0,$ and
\begin{gather*}
\vint_0^\infty\left\langle\vi\frac{\partial u}{\partial t}+\Delta u+Vu+ag(u),\vphi\right\rangle_{X_n^\star,X_n}\psi(t)\,\d t
=\vint_0^\infty\big\langle f(t),\vphi\big\rangle_{X_n^\star,X_n}\psi(t)\,\d t.
\end{gather*}
It follows that $u$ satisfies \eqref{defsol21} in $L^p_\loc(0,\infty;X_n^\star),$ hence in $\Dr^\p\big((0,\infty)\times\Omega\big).$ Moreover, \eqref{thmstrongaH13} and \eqref{thmstrongaH14} come from \eqref{lemthmsaH11}, \eqref{lemthmsaH13}, \eqref{lemthmsaH14}, \eqref{lemthmstrongaH23c} and the lower semicontinuity of the norm. Finally, since $\|U\|_{L^\infty((0,\infty)\times\Omega)}\le1,$ \eqref{thmstrongaH15} comes from \eqref{thmstrongaH12}, \eqref{lemVL2} and \eqref{defsol21}, from which we deduce that $u$ solves~\eqref{defsol21} in $ H^{-1}(\Omega)+L^\infty(\Omega),$ for almost every $t>0.$ Now, assume that $p=\infty.$ It follows that the conclusion of the theorem stands for any $p<\infty.$ Changing the null set $N_0$ by a bigger one, if necessary, \eqref{thmstrongaH16} comes from \eqref{thmstrongaH13}, \eqref{defsol22}, \eqref{lemVL2} and \eqref{defsol21}. Finally, it follows from \eqref{thmstrongaH16} and the embedding $L^\infty(\Omega)\inj Y_n^\star,$ that \eqref{thmstrongaH12} holds for $p=\infty,$ and that $u$ satisfies \eqref{defsol21} in $L^\infty_\loc(0,\infty;X_n^\star).$ The theorem is proved.
\medskip
\end{vproof}

\section*{Acknowledgements}
\addcontentsline{toc}{section}{Acknowledgements}
\baselineskip .5cm
P.~Bégout acknowledges funding from ANR under grant ANR-17-EURE-0010 (Investissements d’Avenir program). He is also very grateful to Prof.~Michel 
Talagrand for the useful discussions he had with him, and for having given him Example~\ref{exanm}. The research of J.~I.~D\'{\i}az was partially supported by the project PID2020-112517GB-I00 of the Spain State Research Agency (AEI).

\baselineskip .4cm

%\bibliographystyle{abbrv}
%\bibliography{Paper18}
%\addcontentsline{toc}{section}{References}

\def\cprime{$^\prime$}

\end{document}